\documentclass[english,10pt,a4paper]{article}
\usepackage[utf8]{inputenc}
\usepackage[T1]{fontenc}
\usepackage{babel}
\usepackage{hyperref}
\usepackage{amsmath}
\usepackage{amssymb}
\usepackage{amsthm}
\usepackage{textcomp}
\theoremstyle{definition}

\newtheorem{Definition}{Definition}[section]

\newtheorem*{Proof}{Proof}
\newtheorem{Proposition}{Proposition}[section]
\newtheorem{Corollary}{Corollary}[section]
\newtheorem{Conjecture}{Conjecture}

\newtheorem{Remark}{Remark}
\newtheorem{Example}{Example}
\usepackage{authblk}
\usepackage[numbers,sort&compress]{natbib}
\usepackage{tocbibind}
\usepackage{tikz}
\usetikzlibrary{decorations.pathreplacing}
\usepackage{tikz-cd}
\usepackage{graphicx}
\usepackage{listings}
\usepackage{xcolor}
\usepackage{caption}
\usepackage{url}
\usepackage{orcidlink}
\numberwithin{figure}{section}
\begin{document}
	
	\title{The chromatic noncommutative symmetric function of oriented trees}
	\author[1]{Lingxiao Hao\thanks{ Corresponding author:  166015167@qq.com}\orcidlink{0009-0009-1242-8956}}	
	\author[2]{Shenglin Zhu\orcidlink{0009-0005-0888-4271}}
	\affil[1]{School of Mathematics and Statistics, Lanzhou University, Lanzhou 730000, China}
	\affil[2]{School of Mathematical Sciences, Fudan University, Shanghai 200433, China} 
	\date{\today}
	\maketitle
	\begin{abstract}
		A long-standing question is whether chromatic symmetric functions can distinguish non-isomorphic trees. Campbell introduced chromatic noncommutative
			symmetric functions for digraphs, which lift chromatic
		symmetric functions to NSym, and asked to what extent they can distinguish non-isomorphic oriented trees. In this article, we prove that chromatic noncommutative
			symmetric functions can reconstruct oriented stars, oriented double stars, some oriented caterpillars, and some oriented paths.  
                      \end{abstract}
                       \text{Keywords:  } chromatic symmetric function;  chromatic noncommutative symmetric function;  digraph; tree; star; double star;  caterpillar; oriented path
	\section{Introduction}
	In 1995, Stanley~\cite{STANLEY1995166} introduced \emph{chromatic symmetric functions} in \text{Sym} for undirected graphs and raised a well-known conjecture. 
	\begin{Conjecture}\label{CSF}
		Non-isomorphic trees have different chromatic symmetric functions.
	\end{Conjecture}
	Many significant results concerning this conjecture have emerged. It has been proved for paths, stars, squids, and spiders by J. L. Martin, M. Morin, and J. D. Wagner~\cite{MARTIN2008237}, for caterpillars by M. Loebl and J-S. Sereni~\cite{AIHPD_2019__6_3_357_0}, for trees with up to 29 vertices by S. Heil and C. Ji~\cite{heil2018algorithmcomparingchromaticsymmetric}, for proper trees of diameter at most 5 by J. Aliste-Prieto et al.~\cite{doi: 10.1137/22M148046X}, and for trees with exactly two vertices of degree at least 3 by Y. Wang et al.~\cite{wang2024classtreesdeterminedchromatic}.

In 2024, Campbell~\cite{campbell2024liftchromaticsymmetricfunctions} introduced chromatic noncommutative
			symmetric functions by lifting chromatic
	symmetric functions to \text{NSym}, and asked to what extent these functions can distinguish non-isomorphic digraphs whose underlying graphs are trees (referred to as oriented trees).

 We prove that chromatic noncommutative
			symmetric functions can reconstruct oriented stars, oriented double stars, some oriented caterpillars, and some oriented paths. We also apply known positive results on Conjecture~\ref{CSF} to determine the underlying graph of these digraphs, and determine the directions of edges---from leaf edges to non-leaf ones---based on the coefficients of chromatic noncommutative
			symmetric functions.

	The paper is organized as follows.

                        In Section~\ref{star double star and caterpillar}, we prove that all oriented stars, oriented double stars, and some oriented caterpillars can be reconstructed from their chromatic noncommutative
			symmetric functions.

	In Section~\ref{path}, we present some results on oriented paths. In detail, we prove that chromatic noncommutative symmetric functions can reconstruct directed paths, oriented paths that are almost symmetric, unimodal, or inflection-dense, and oriented paths with up to 24 vertices. For an oriented path $D$, we focus on compositions where exactly one part is larger than 1 in $\textbf{X}_{D}$, and show that they (or rather, their coefficients) can reveal substantial information about $D$. Based on this observation, we present an algorithm in Python (shown in Appendix~\ref{app sec A}) and prove that chromatic noncommutative symmetric functions can distinguish oriented paths with up to 24 (and possibly more) vertices.

	\section{Preliminaries}\label{Preliminaries}
	A \emph{composition} is a finite sequence of positive integers $\alpha=(\alpha_{1},\alpha_{2},\ldots,\alpha_{l(\alpha)})$, where $l(\alpha)$ is called its \emph{length} and $\alpha_{1},\alpha_{2},\ldots,\alpha_{l(\alpha)}$ are called its \emph{parts}. A composition $\alpha$ is called a \emph{rearrangement} of another composition $\beta$ if the multiset of parts of $\alpha$ equals the multiset of parts of $\beta$. A \emph{partition} is a composition whose parts are non-increasing.

	The ring of \emph{symmetric functions} $\text{Sym}$~\cite{10.1093/oso/9780198534891.001.0001} consists of polynomials in commuting variables (i.e., $x_{i} x_{j} = x_{j} x_{i}$) that are invariant under any permutation of the variables. Formally, a polynomial $f \in \mathbb{Z}[x_{1}, x_{2}, \dots]$ is in $\text{Sym}$ if $f(x_{\sigma(1)}, x_{\sigma(2)}, \dots) = f(x_{1}, x_{2}, \dots)$ for any permutation $\sigma$ of indices. In this article, we need to introduce the \emph{power sum basis} $\left\lbrace p_{\lambda}\right\rbrace $ of $\text{Sym}$. Let $p_{n}= \sum_{i \geq 1} x_{i}^{n}$. For a partition $\lambda=(\lambda_{1},\lambda_{2},\ldots,\lambda_{l})$, $p_{\lambda}=p_{\lambda_{1}}p_{\lambda_{2}}\cdots p_{\lambda_{l}}.$

	The algebra of \emph{noncommutative symmetric functions} $\operatorname{NSym}$ is the free associative algebra $\mathbb{K}\langle H_1, H_2, \dots\rangle$ generated by the corresponding elements $H_1, H_2,\\ \dots$\cite{gelfand1994noncommutativesymmetricfunctions}.

The \emph{chromatic symmetric function} (\emph{CSF}) is a graph-associated symmetric function that encodes coloring information of a graph. For a simple graph $G = (V, E)$ with vertices $V = \left\lbrace v_{1}, v_{2}, \dots, v_{n}\right\rbrace$ and edges $E$, its CSF~\cite[Definition 2.1]{STANLEY1995166} is defined over commuting variables $X = \left\lbrace x_{1}, x_{2}, \dots\right\rbrace$. 
	$X_{G} = \sum_{\kappa:  V \to \mathbb{N}} \prod_{v \in V} x_{\kappa(v)}$
	where the sum runs over all proper colorings $\kappa$ of $G$ (i.e., $\kappa(u) \neq \kappa(v)$ for every edge $(u, v) \in E$), and $\prod_{v \in V} x_{\kappa(v)}$ is the product of variables indexed by the color assigned to each vertex.
	The expansion of the CSF in the power sum basis is as follows (see~\cite[Theorem 2.5]{STANLEY1995166}).
	$$X_{G}=\mathop{\sum}\limits_{S\subseteq E(G)}(-1)^{|S|}p_{\lambda(S)}.$$
	For each subset $S\subseteq E(G)$, the parts of the partition $\lambda(S)$ are the (vertex) sizes of the connected components of the spanning subgraph of $G$ with edges $S$.

        CSF is defined for undirected graphs. For \emph{directed graphs}, also termed \emph{digraphs}, an edge $(v,w)$ in a digraph is said to be \emph{from $v$ to $w$}, denoted by
  \begin{center}
    \begin{tikzcd}
      v \arrow[r] & w.
    \end{tikzcd}
    \captionof{figure}{}
  \end{center}

        Throughout this paper, we say in a digraph $D$ that a vertex $v$ is \emph{joined with} a vertex $w$ if $v$ is adjacent to $w$ in the underlying undirected graph of $D$. Two digraphs $D$ and $D'$ are said to be \emph{isomorphic} if $D'$ can be obtained by relabeling the vertices of $D$. We regard isomorphic digraphs as identical. For a vertex $v$ in a digraph, the \emph{outdegree} of $v$ is the number of edges from $v$, and the \emph{indegree} of $v$ is the number of edges to $v$.

Let $D=\left(V(D), E(D)\right)$ be a digraph. If $S$ is a subset of $V(D)$, the \emph{total degree} of $S$ is defined in~\cite[Definition 1]{campbell2024liftchromaticsymmetricfunctions} as $$\text{td}S=\mathop{\sum}\limits_{v\in S}\text{od}v-\mathop{\sum}\limits_{v\in S}\text{id}v, $$ where $\text{od}$ and $\text{id}$ denote the outdegree and indegree. For convenience, we abbreviate $\text{td}\left\lbrace v\right\rbrace$ as $\text{td}v$ for a single vertex $v$. A vertex $v$ in a digraph is a \emph{source} if $\text{td}v=\text{od}v$, and is a \emph{sink} if $\text{td}v=-\text{id}v$.

If $S$ is a subset of $E(D)$, the composition $\alpha(S)$ of $S$ is defined in~\cite[Definition 3]{campbell2024liftchromaticsymmetricfunctions} as $(|C^{S}_{1}|,\ldots,|C^{S}_{l}|)$, where $C^{S}_{1},\ldots,C^{S}_{l}$ are the connected components of the spanning subgraph of $D$ induced by $S$, and $i<j$ if $\text{td}C^{S}_{i}>\text{td}C^{S}_{j}$, or $\text{td}C^{S}_{i}=\text{td}C^{S}_{j}$ and $|C^{S}_{i}|>|C^{S}_{j}|$. For this notion, one can check the following example.
	\begin{Example}
		Let $D$ be the digraph
\begin{center}
  \begin{tikzcd}
    v_{1} \arrow[r, "e_{1}"] \arrow[d, "e_{4}"] & v_{2} \arrow[d, "e_{2}"] &    \\
    v_{3} \arrow[r, "e_{3}"]                 & v_{4} \arrow[r, "e_{5}"] & v_{5},
  \end{tikzcd}
  \captionof{figure}{}
\end{center}
	and $S=\left\lbrace e_{2}, e_{3}\right\rbrace $. There are three connected components $\left\{v_{1}\right\}$, $\left\{v_{2},v_{3},v_{4}\right\}$ and $\left\{v_{5}\right\}$, with total degrees $\text{td}v_{1}=+2$, $\text{td}\left\lbrace v_{2},v_{3},v_{4}\right\rbrace =-1$ and $\text{td}v_{5}=-1$, respectively. Thus, $\alpha(S)=(|\left\lbrace v_{1}\right\rbrace |, |\left\lbrace v_{2},v_{3},v_{4}\right\rbrace|, |\left\lbrace v_{5}\right\rbrace |)=(1,3,1)$.
	\end{Example}
As shown in the example, for $S\subseteq E(D)$ we assign to $\alpha=\alpha(S)=(|C^{S}_{1}|,\ldots,|C^{S}_{l}|)$ a \emph{degree} defined as $(\text{td}C^{S}_{1},\ldots,\text{td}C^{S}_{l})$. However, to obtain $\alpha(S)$ we need to compute $(\text{td}C^{S}_{1},\ldots,\text{td}C^{S}_{l})$ first.  Note that one composition $\alpha$ may have several degrees, since the same $\alpha$ may arise from different $S$. 
	\begin{Definition}[\protect{\cite[Definition 4]{campbell2024liftchromaticsymmetricfunctions}}]
		Let $D$ be a digraph. The \emph{chromatic noncommutative
			symmetric function} (referred to as \emph{CNSF}) of $D$ is defined as\\
		$$\textbf{X}_{D}=\mathop{\sum}\limits_{S\subseteq E(D)}(-1)^{|S|}\Psi_{\alpha(S)}.$$
	\end{Definition}
For a composition $\alpha=(\alpha_{1},\alpha_{2},\ldots,\alpha_{l})$, $\Psi_{\alpha}$ is defined as $\Psi_{\alpha}=\Psi_{\alpha_{1}}\Psi_{\alpha_{2}}\cdots\Psi_{\alpha_{l}}$, where $\Psi_{n}$ is referred to as a power sum of the first kind (for its definition, see~\cite[Definition 2.1]{duchamp: hal-00018540}). Note that the forgetful map $\chi: \text{NSym}\rightarrow \text{Sym}$ sends $\Psi_{n}$ to $p_{n}$; hence the following proposition holds.
	\begin{Proposition}[\protect{\cite[Theorem 2]{campbell2024liftchromaticsymmetricfunctions}}]\label{1}
		For a digraph $D$ with $G$ as its underlying graph, the
		projection identity $\chi(\textbf{X}_{D}) = X_{G}$ holds.
	\end{Proposition} 
The above proposition allows us to apply the results about Conjecture~\ref{CSF} to determine the underlying graph of a digraph $D$ when we only know $\textbf{X}_{D}$.  When the underlying graph of $D$ is a tree with $n$ vertices, we find that $$|S|=n-l(\alpha(S))$$ for any subset $S$ of $V(D)$. Therefore, when $\alpha(S)=\alpha(S')$, the terms corresponding to $S$ and $S'$ do not cancel. In this case, we rewrite $\textbf{X}_{D}$ as $$\textbf{X}_{D}=\mathop{\sum}(-1)^{n-l(\alpha)}c(\alpha)\Psi_{\alpha},$$ where $c(\alpha)$ is the number of $S$ such that $\alpha(S)=\alpha$. For brevity, we may write $c((\alpha_{1},\ldots,\alpha_{l}))$ as $c(\alpha_{1},\ldots,\alpha_{l})$. We say a composition $\alpha$ \emph{appears} in $\textbf{X}_{D}$ if $c(\alpha)\neq 0$.
	\begin{Remark}
Note that the degree of any composition in $\textbf{X}_{D}$ is a weakly decreasing sequence. If the underlying graph of $D$ is a tree, then for any composition $\alpha$ in $\textbf{X}_{D}$ such that $l(\alpha)\geq 2$, the first part of the degree of $\alpha$ must be positive, and the last part of the degree must be negative. To see this, note that by contracting each connected component of the spanning subgraph of $D$ induced by an edge subset $S$ into a single vertex, we obtain a non-trivial tree $D_{S}$. Each vertex of $D_{S}$ retains the total degree of the corresponding connected component. There must be some leaves in $D_{S}$, and their total degrees must be non-zero.
	\end{Remark}
An edge is called a \emph{leaf edge} if it is incident to a leaf, and it is called a \emph{non-leaf edge} otherwise. To distinguish the directions of leaf edges, we give the following definition.
	\begin{Definition}
          Let $v$ be a leaf of a tree. We call $v$ \emph{positive} if $\text{td}v =+1$, and \emph{negative} if $\text{td}v=-1$.
	\end{Definition}
	
	\section{Star, Double Star and Caterpillar}\label{star double star and caterpillar}
In this section we study how digraphs whose underlying graphs are stars, double stars, and caterpillars (hereafter called oriented stars, oriented double stars, and oriented caterpillars, for convenience) can be distinguished by their CNSFs.

The \emph{star} $S_{n}$ is a tree consisting of $n-1$ leaves and a unique non-leaf vertex, which is called the \emph{center}. The \emph{double star} $S_{a,b}$ is a tree with exactly two non-leaf vertices: one is adjacent to $a$ leaves, and the other to $b$ leaves.

We can obtain the number of positive and negative leaves of an oriented tree from its CNSF.
	
	\begin{Proposition}
         Let $D$ be an oriented tree with $n$ vertices. There are exactly $c(n-1,1)$ negative leaves and $c(1,n-1)$ positive ones in $D$.
	\end{Proposition}
	\begin{Proof}
Every $\alpha$ of length 2 appearing in $\textbf{X}_{D}$ is obtained by deleting one edge, since $D$ is an oriented tree. If one part of $\alpha$ is 1, then the deleted edge must be a leaf edge. The degree of $\alpha$ must be $(+1,-1)$.

		If $n>2$, we can determine whether the leaf is positive or negative from the position of 1 in $\alpha$.

		If $n=2$, then there is one positive leaf and one negative leaf, and $c(1,1)=1.$
		$\qedsymbol$
	\end{Proof}
From this proposition, we know the number of all leaves of $D$ from $\textbf{X}_{D}$. If an oriented tree $D$ with $n$ vertices has exactly $n-1$ leaves (respectively $n-2$ leaves), then $D$ must be an oriented star (respectively oriented double star). Moreover, we can reconstruct an oriented star from the number of positive and negative leaves.
	\begin{Corollary}
		An oriented star $D$  can be reconstructed from its $\textbf{X}_{D}$.
	\end{Corollary}
	For oriented double stars, the situation is much more complicated. Suppose the unique non-leaf edge is from $v_{0}$ to $v_{1}$. For each $i=0,1$, denote the number of positive leaves and negative leaves joined to $v_{i}$ by $n_{+}(v_{i})$ and $n_{-}(v_{i})$.
	\begin{Proposition}\label{double star}
		Let $D$ be an oriented double star. Then it can be reconstructed from $\textbf{X}_{D}$.
	\end{Proposition}
	\begin{Proof}
There is only one composition of length two with both parts larger than 1 in $\textbf{X}_{D}$ since there is only one non-leaf edge in $D$. Let this composition be $(a+1,b+1)$. The unique non-leaf edge is from the center $v_{0}$ of $S_{a+1}$ to the center $v_{1}$ of $S_{b+1}$, as follows: \\
                \begin{center}
                    \begin{tikzcd}
                      {} \arrow[rd] & \cdots                     & \cdots           & {} \arrow[ld] \\
                      \cdots        & v_{0} \arrow[r] \arrow[ld] & v_{1} \arrow[rd] & \cdots        \\
                      {}            & \cdots                     & \cdots           & {}           
                    \end{tikzcd}
                    \captionof{figure}{}
                \end{center}

                Consider the compositions 			
                obtained by deleting the non-leaf edge and a leaf edge.

                There are 4 situations:
                \begin{enumerate}
						\item If the leaf edge is from $v_{0}$, then the corresponding composition is $(a,b+1,1)$ with degree $(+2, -1, -1)$.
						\item If the leaf edge is to $v_{0}$, then the composition is $(1,a,b+1)$ with degree $(+1,0, -1)$.
						\item If the leaf edge is from $v_{1}$, then the composition is $(a+1,b,1)$ with degree $(+1,0, -1)$.
						\item If the leaf edge is to $v_{1}$, then the composition is $(a+1,1,b)$ with degree $(+1, +1, -2)$.
                                                \end{enumerate}

                \item  Case A:  $a,b\geq 2$
		\begin{quote}
The 4 compositions are distinct even if $a=b$. In addition, compositions of length 3 with two parts larger than 1 in $\textbf{X}_{D}$ can only be obtained by deleting the non-leaf edge and a leaf edge. Hence, $n_{-}(v_{0})=c(a,b+1,1)$, $n_{+}(v_{0})=c(1,a,b+1)$, $n_{-}(v_{1})=c(a+1,b,1)$, and $n_{+}(v_{1})=c(a+1,1,b)$.
		\end{quote}	
		\item Case B:  $a=1$, $b\geq 2$ 
		\begin{quote}
The composition $(2,b,1)$ can only be obtained by deleting the non-leaf edge and a leaf edge from $v_{1}$. There are $3+b$ vertices in total, and there are $c(b+2,1)$ negative leaves. Hence, $n_{-}(v_{1})=c(2,b,1)$, $n_{+}(v_{1})=b-n_{-}(v_{1})$, $n_{-}(v_{0})=c(b+2,1)-n_{-}(v_{1})$, and $n_{+}(v_{0})=1-n_{-}(v_{0})$.
		\end{quote}	
		\item Case C:  $b=1$, $a\geq 2$ \begin{quote}
The composition $(a,2,1)$ can only be obtained by deleting the non-leaf edge and a leaf edge from $v_{0}$. There are $3+a$ vertices in total and $c(a+2,1)$ negative leaves. Hence, $n_{-}(v_{0})=c(a,2,1)$, $n_{+}(v_{0})=a-n_{-}(v_{0})$, $n_{-}(v_{1})=c(a+2,1)-n_{-}(v_{0})$, and $n_{+}(v_{1})=1-n_{-}(v_{1})$.
		\end{quote}	
		\item Case D:  $a = b = 1$
		\begin{quote}
In this case, there are 4 possibilities. Denote the leaf joined with $v_{0}$ by $v'_{0}$ and the leaf joined with $v_{1}$ by $v'_{1}$.
			
			We further distinguish three mutually exclusive subcases:
			\begin{enumerate}
				\item Subcase D.1:  If $c(3,1) = 2$, then both leaves are negative.
				\item Subcase D.2:  If $c(1,3) = 2$, then both leaves are positive.
				\item Subcase D.3:  If $c(3,1) = c(1,3) = 1$, then one leaf is positive while the other is negative.
				
				Further, under Subcase D.3, the composition $(1,1,2)$ can only be obtained by deleting the non-leaf edge and a leaf edge incident to $v_{0}$: 
				\begin{itemize}
					\item  If $c(1,1,2) = 1$, then $v'_{0}$ is positive while $v'_{1}$ is negative.
					\item   If $c(1,1,2) = 0$, then $v'_{1}$ is positive while $v'_{0}$ is negative.
				\end{itemize}
			\end{enumerate}
		\end{quote}	
$\qedsymbol$
	\end{Proof}

A \emph{caterpillar} is a tree in which every vertex not on a single path is a leaf. This path is called the \emph{spine}. For a caterpillar $T$ with spine vertices $v_{0},v_{1},\ldots,v_{l}$, the \emph{leaf number} $f_{i}$ of $v_{i}$ is the number of leaves adjacent to $v_{i}$ (see~\cite{MARTIN2008237}).
	Before discussing general oriented caterpillars, we first present a result demonstrating that CSFs can distinguish certain caterpillars.
	\begin{Proposition}[\protect{\cite[Theorem 11]{MARTIN2008237}}]\label{caterpillar}
	Let $T$ be a caterpillar with strictly positive and distinct leaf numbers $f_i$. Then $T$ can be uniquely reconstructed from its characteristic subtree function $X_T$.
	\end{Proposition}
	Let $a_{i}=f_{i}+1$ and $a_{j,k}=\mathop{\sum}\limits_{j \leq i\leq k}a_{i}$.	With the above proposition, we obtain the following result.
	\begin{Proposition}
		Let $D$ be an oriented caterpillar of spine length $l$ where all the numbers $f_{0},\ldots,f_{l}$ and the $a_{j,k}$ ($0\le j\le k\le l$) of its underlying caterpillar are pairwise distinct and larger than 1. Then $D$ can be reconstructed from $\textbf{X}_{D}$.
              \end{Proposition}
        
	\begin{Proof}
According to Propositions~\ref{1} and~\ref{caterpillar}, we can obtain the underlying caterpillar of $D$. Denote its spine vertices by $v_{0},\ldots, v_{l}$. Now we know all $f_{i}$ and $a_{j,k}$. Denote the edge between $v_{i-1}$ and $v_{i}$ by $e_{i}$. The conditions in this proposition ensure that every composition discussed below corresponds to only one way of deleting edges. 

We first discuss the non-leaf edges. For each $0\leq i \leq l-1$, find the unique $\alpha_{i}$ in $\textbf{X}_{D}$ that is a rearrangement of $(a_{0,i},a_{i+1,l})$. The composition $\alpha_{i}$ must be obtained by deleting $e_{i+1}$, and the degree of $\alpha_{i}$ must be $(+1,-1)$. If $\alpha_{i}=(a_{0,i},a_{i+1,l})$, then $e_{i+1} $ is from $v_{i}$ to $v_{i+1}$; conversely, if $\alpha_{i}=(a_{i+1,l},a_{0,i})$, then $e_{i+1} $ is from $v_{i+1}$ to $v_{i}$.\\

Now we need to determine the directions of the leaf edges. For each $i$, denote the number of positive leaves and negative leaves joined with $v_{i}$ by $n_{+}(v_{i})$ and $n_{-}(v_{i})$ --- we only need to know one of them since $n_{+}(v_{i})+n_{-}(v_{i})=f_{i}$.
		
For a leaf attached to $v_0$, deleting its incident edge and $e_1$ yields three components with $1$, $f_0$, and $a_{1,l}$ vertices, respectively.
\begin{enumerate}
	\item If the leaf is positive, the total degree of the component with $f_0$ vertices is either $0$ or $-2$, as shown below.
\\
			{
\centering
\begin{tikzcd}
				{} \arrow[rd] & \cdots                   &        & {} \arrow[rd] & \cdots &                             \\
				& v_{0} \arrow[r, "e_{1}"] & \cdots & \text{or}            & v_{0}  & \cdots. \arrow[l, "e_{1}"']
			\end{tikzcd}
\captionof{figure}{}
}The corresponding composition is either $(1,f_{0},a_{1,l})$ with degree $(+1, 0, -1)$, or $(a_{1,l},1,f_{0})$ with degree $(+1, +1, -2)$. \item If it is a negative leaf, then the total degree of the component with $f_{0}$ vertices is either $0$ or $+2$, as shown below:\\
			{
\centering
\begin{tikzcd}
				{} & \cdots           &                            & {} & \cdots                              &         \\
				& v_{0} \arrow[lu] & \cdots \arrow[l, "e_{1}"'] & \text{or} & v_{0} \arrow[r, "e_{1}"] \arrow[lu] & \cdots.
			\end{tikzcd}
\captionof{figure}{}
}  The corresponding composition is either $(a_{1,l},f_{0},1)$ with degree $(+1, 0, -1)$ or $(f_{0},a_{1,l},1)$ with degree $(+2, -1, -1)$. 
		\end{enumerate} 

We claim that a composition in $\textbf{X}_{D}$ whose multiset of parts is $\{ f_{0}, a_{1,l}, 1\}$ must be obtained by deleting a leaf edge incident to $v_{0}$ and the edge $e_{1}$. From this claim we obtain
\[
n_{+}(v_{0}) = c(1, f_{0}, a_{1,l}) + c(a_{1,l}, 1, f_{0}),
\]
and similarly
\[
n_{+}(v_{l}) = c(1, f_{l}, a_{0,l-1}) + c(a_{1,l}, 1, a_{0,l-1}).
\]

We now prove the claim. If only non-leaf edges are deleted, each component size must be of the form $a_{j,k}$ for some $j\le k$. Consider a composition in $\textbf{X}_{D}$ with multiset $\{ f_{0}, a_{1,l}, 1\}$. The isolated vertex must be a leaf, because isolating a non-leaf vertex $v_i$ would produce $a_i$ singletons. The composition is obtained by deleting one leaf edge and one non-leaf edge. If we delete only that non-leaf edge, the component sizes would be $f_{0}+1 = a_{0}=a_{0,0}$ and $a_{1,l}$; they cannot be $f_{0}$ and $a_{1,l}+1$. Hence the non-leaf edge must be $e_{1}$, and the leaf edge is incident to $v_{0}$.

A similar argument shows that for $1\le i\le l-1$, a composition in $\textbf{X}_{D}$ whose multiset of parts is $\{ f_{i}, a_{0,i-1}, a_{i+1,l}, 1\}$ must be obtained by deleting $e_{i}$, $e_{i+1}$, and a leaf edge incident to $v_{i}$.

For a leaf attached to $v_i$ with $1\le i\le l-1$, deleting its incident edge, $e_i$, and $e_{i+1}$ yields four components with $f_i$, $a_{0,i-1}$, $a_{i+1,l}$, and $1$ vertices, respectively. We consider three cases depending on $e_i$ and $e_{i+1}$:
\begin{enumerate}\item  Case A. Both $e_{i+1}$ and $e_{i}$ are from $v_{i}$.\begin{quote}
						The digraph can be seen as follows: \\
								{
\centering
\begin{tikzcd}
									& {}      &                                                           & {} \arrow[ld] &        \\
									\cdots & v_{i-1} & v_{i} \arrow[l, "e_{i}"'] \arrow[r, "e_{i+1}"] \arrow[lu] & v_{i+1}       & \cdots. \\
									&         & \cdots                                                    &               &       
								\end{tikzcd}
\captionof{figure}{}
}\begin{enumerate}
\item The composition obtained by deleting $e_{i+1}$, $e_{i}$ and a leaf edge from $v_{i}$ is either $(f_{i},a_{0,i-1},a_{i+1,l},1)$ or $(f_{i},a_{i+1,l},a_{0,i-1},1)$ (depending on whether $a_{0,i-1}>a_{i+1,l}$ or $a_{0,i-1}<a_{i+1,l}$) with degree $(+3, -1, -1,-1)$.
\item The composition obtained by deleting $e_{i+1}$, $e_{i}$ and a leaf edge to $v_{i}$ is either $(f_{i},1,a_{0,i-1},a_{i+1,l})$ or $(f_{i},1,a_{i+1,l},a_{0,i-1})$ (depending on whether $a_{0,i-1}>a_{i+1,l}$ or $a_{0,i-1}<a_{i+1,l}$) with degree $(+1, +1, -1,-1)$.
\end{enumerate} Hence, $n_{+}(v_{i})=c(f_{i},1,a_{0,i-1},a_{i+1,l})+c(f_{i},1,a_{i+1,l},a_{0,i-1})$ in case A.
						\end{quote}
	\item 	Case  B. Both $e_{i+1}$ and $e_{i}$ are to $v_{i}$.\begin{quote}
				The digraph can be displayed as follows: \\
						{
\centering
\begin{tikzcd}
							& {}                           &                  & {} \arrow[ld]               &        \\
							\cdots & v_{i-1} \arrow[r, "e_{i}"] & v_{i} \arrow[lu] & v_{i+1} \arrow[l, "e_{i+1}"'] & \cdots. \\
							&                              & \cdots           &                             &       
						\end{tikzcd}
\captionof{figure}{}
}
						\begin{enumerate}
							\item The composition obtained by deleting $e_{i+1},e_{i}$ and a leaf edge from $v_{i}$ is either $(a_{0,i-1},a_{i+1,l},f_{i},1)$ or $(a_{i+1,l},a_{0,i-1},f_{i},1)$ (depending on whether $a_{0,i-1}>a_{i+1,l}$ or $a_{0,i-1}<a_{i+1,l}$) with degree $(+1, +1, -1,-1)$.
							\item The composition obtained by deleting $e_{i+1}$, $e_{i}$ and a leaf edge to $v_{i}$ is either $(a_{0,i-1},a_{i+1,l},1,f_{i})$ or $(a_{i+1,l},a_{0,i-1},1,f_{i})$ (depending on whether $a_{0,i-1}>a_{i+1,l}$ or $a_{0,i-1}<a_{i+1,l}$) with degree $(+1, +1, +1,-3)$. \end{enumerate}Hence, $n_{+}(v_{i})=c(a_{0,i-1},a_{i+1,l},1,f_{i})+c(a_{i+1,l},a_{0,i-1},1,f_{i})$ in case B.
				\end{quote} 
		\item Case C. One of $e_{i+1}$ and $e_{i}$ is to $v_{i}$ while the other is from $v_{i}$. \begin{quote}
		The digraph is presented below: \\
				{
\centering
\begin{tikzcd}
					& {}                         &                                         & {} \arrow[ld] &        \\
					\cdots & {} \arrow[r, "\text{non-leaf} "'] & v_{i} \arrow[r, "\text{non-leaf}"'] \arrow[lu] & {}            & \cdots. \\
					&                            & \cdots                                  &               &       
				\end{tikzcd}
\captionof{figure}{}
} Deleting $e_{i+1}$, $e_{i}$ and a leaf edge incident to $v_{i}$ yields a composition $\alpha$ whose multiset of parts is $\left\lbrace a_{0,i-1},a_{i+1,l},1,f_{i}\right\rbrace $ and whose degree is $(+1, +1, -1,-1)$. We can determine the total degree of the leaf by observing the position of $1$ in $\alpha$.
Hence, in case C, $$n_{+}(v_{i})=\mathop{\sum}\limits_{\left\lbrace a,b,c\right\rbrace =\left\lbrace a_{0,i-1},a_{i+1,l},f_{i}\right\rbrace}c(a,1,b,c).$$
		\end{quote}\end{enumerate}
		$\qedsymbol$
	\end{Proof}

\section{Path}\label{path}
In this section, we consider special oriented trees whose underlying graphs are paths; we call them \emph{oriented paths}. Throughout this section, let $D$ be an oriented path with $n$ vertices. The vertices are denoted $v_{0},v_{1},\ldots,v_{n-1}$ (in some proofs we may designate a specific leaf as $v_{0}$), and the edge between $v_{i-1}$ and $v_{i}$ is denoted $e_{i}$ for $1\leq i\leq n-1$. For $0 \le i < j \le n-1$, let $p_{i,j}$ denote the digraph induced by $v_{i},v_{i+1},\ldots,v_{j}$. Then, for $1 \le i \le j \le n-2$, we have $\operatorname{td} p_{0,i} \in \{-1,+1\}$, $\operatorname{td} p_{j,n-1} \in \{-1,+1\}$, and $\operatorname{td} p_{i,j} \in \{-2,0,+2\}$. We first introduce a proposition that can be used to determine whether the underlying graph of a digraph is a path.
	\begin{Proposition}[\protect{\cite[Corollary 5]{MARTIN2008237}}]\label{degree and path sequence}
		The degree and path sequences of a tree $T$ can be recovered from its CSF $X_{T}$.
	\end{Proposition}

\subsection{Directionality}\label{Directionality}

 We first consider the most restricted kind of oriented paths: directed paths. A \emph{directed path} is an oriented path $D$ consisting of vertices $v_0, v_1, \dots, v_k$ and edges $e_i$ from $v_{i-1}$ to $v_i$ for all $1 \leq i < k$. The \emph{length} of $D$ is $k$, and we say $D$ goes from $v_0$ to $v_k$.  
	{
\centering
\begin{tikzcd}
		v_{0} \arrow[r, dashed] & v_{k}.
	\end{tikzcd}
\captionof{figure}{}
}

	\begin{Proposition}\label{directed path}
		If $D$ is a directed path, then $D$ can be uniquely determined by $\textbf{X}_{D}$.
		More precisely, if $D$ is an oriented path with $n$ vertices, then $D$ is a directed path if and only if, for each $k\leq \frac{n}{2}$, both $(1^{k},n-k)$ and $(n-k,1^{k})$ appear in $\textbf{X}_{D}$.
	\end{Proposition}
	\begin{Proof}
		
		If $D$ is a directed path with $n$ vertices, then both $(1,n-1)$ and $(n-1,1)$ appear in  $\textbf{X}_{D}$. From the unique source to the unique sink, denote the vertices by $v_{0},v_{1},\ldots,v_{n-1}$ and the edges by $e_{1},\ldots,e_{n-1}$. For $k\leq \frac{n}{2}$, deleting $e_{1},e_{2},\ldots,e_{k}$ yields $(1^{k},n-k)$ with degree $(+1,0,0,\ldots, 0,-1)$, and deleting $e_{n-k},\ldots,e_{n-1}$ yields $(n-k,1^{k})$ with degree $(+1,0,0,\ldots, 0,-1)$. \\
		
		Now we prove the other direction.  Observe the following diagram: \\
		{
\centering
\begin{tikzcd}
			v_{0} \arrow[rr, dashed] &  & v_{k} \arrow[r] & v_{k+1} & v_{n-k-2} \arrow[r] & v_{n-k-1} \arrow[rr, dashed] &  & v_{n-1}.
		\end{tikzcd}
\captionof{figure}{}
} That both $(1,n-1)$ and $(n-1,1)$ appear in $\textbf{X}_{D}$ tells us that one leaf of $D$ is positive and the other is negative. From the positive leaf to the negative one, denote the vertices by $v_{0},v_{1},\ldots,v_{n-1}$ and the edges by $e_{1},\ldots,e_{n-1}$.
		Now we have that $e_{1}$ is from $v_{0}$ to $v_{1}$ and $e_{n-1}$ is from $v_{n-2}$ to $v_{n-1}$. If $n\leq 3$, then $D$ is a directed path. Now consider the case $n\geq 4$. \begin{quote}
		
		For $1\leq k \leq\frac{n}{2}-1$, $n-k-1\geq \frac{n}{2}-2> 1$.
		Suppose that for each $1\leq i\leq k$,  $e_{i}$ is from $v_{i-1}$ to $v_{i}$ and $e_{n-i}$ is from $v_{n-i-1}$ to $v_{n-i}$, then $\text{td} p_{0,n-k-2},\text{td} p_{k+1,n-1}\in \left\lbrace  +1,-1\right\rbrace  $, and $\text{td} p_{j,j+n-k-2}=0$ for $1\leq j\leq k$.\begin{enumerate}
		\item 
		 Consider the composition $(1^{k+1},n-k-1)$. If one of the isolated vertices is $v_{n-1}$, whose total degree is $-1$, then the total degree of the component of $n-k-1$ vertices must be $-2$, which is a contradiction. Thus, $(1^{k+1},n-k-1)$ must be obtained by deleting $e_{1},\ldots,e_{k+1}$ --- the isolated vertices are $v_{0},\ldots,v_{k}$. So $\text{td}p_{k+1,n-1}=-1$, which means $e_{k+1}$ is from $v_{k}$ to $v_{k+1}$.\item
		Now consider $(n-k-1,1^{k+1})$. The total degree of the component of $n-k-1$ vertices must be positive. It can only be $p_{0,n-k-2}$ since $\text{td}p_{k+1,n-1}=-1$ and $\text{td} p_{j,j+n-k-2}=0$ for $j\leq 1\leq k$. Thus, $\text{td}p_{0,n-k-2}=+1$, which means $e_{n-k-1}$ is from $v_{n-k-2}$ to $v_{n-k-1}$.\end{enumerate}\end{quote}
           If $n=2j+1$, then $p_{0,j}$ is a directed path from $v_{0}$ to $v_{j}$ and $p_{j,2j}$ from $v_{j}$ to $v_{2j}$. If $n=2j$, then $p_{0,j}$ is a directed path from $v_{0}$ to $v_{j}$ and $p_{j-1,2j-1}$ from $v_{j-1}$ to $v_{2j-1}$. Therefore, $D$ is a directed path. $\qedsymbol$
	\end{Proof}
	The proof of the above proposition suggests the following: 
\begin{Proposition}\label{right right}
	Let $D$ be an oriented path with $n$ vertices where $p_{0,k}$ is a directed path from $v_{0}$ to $v_{k}$, and $p_{n-k-1,n-1}$ is a directed path from $v_{n-k-1}$ to $v_{n-1}$ for some $1\leq k<n-2$. Then we can determine the direction of $e_{k+1}$ and $e_{n-k-1}$ from $\textbf{X}_{D}$.
	\end{Proposition}
	\begin{Proof}
Under the condition of this proposition, $\text{td} p_{j,j+n-k-2}=0$ for $j\leq 1\leq k$. We discuss this depending on whether $(1^{k+1},n-k-1)$ appears in $\textbf{X}_{D}$.
		\begin{enumerate}
		\item Case A: The composition $(1^{k+1},n-k-1)$ appears in $\textbf{X}_{D}$. \begin{quote}
		From the proof of Proposition~\ref{directed path}, we have that $e_{k+1}$ is from $v_{k}$ to $v_{k+1}$. In this case, $(n-k-1,1^{k+1})$ appears in $\textbf{X}_{D}$ if and only if $e_{n-k-1}$ is from $v_{n-k-2}$ to $v_{n-k-1}$.
		\end{quote}
		\item Case B: The composition $(1^{k+1},n-k-1)$ does not appear in $\textbf{X}_{D}$. \begin{quote}
		In this case, $e_{k+1}$ is from $v_{k+1}$ to $v_{k}$, as illustrated below: \\
				{
\centering
\begin{tikzcd}[row sep=0.5cm, column sep=0.6cm]
					v_{0} \arrow[r, dashed] & v_{k} & v_{k+1} \arrow[l] & v_{n-k-2} & v_{n-k-1} \arrow[r, dashed] & v_{n-1}.
				\end{tikzcd}
\captionof{figure}{}
}
				In this case, $\text{td}p_{k+1,n-1}=+1$. Deleting $e_{1}\cdots e_{k+1}$ yields $(n-k-1,1^{k+1})$ with degree $(+1,+1,0,\ldots,0,-2)$. \begin{enumerate}
						\item Subcase B.1 If $e_{n-k-1}$ is from $v_{n-k-2}$ to $v_{n-k-1}$, then $\text{td}p_{0,n-k-2}=+1$. In this case, deleting $e_{n-k-1},\cdots, e_{n-1}$ also yields $(n-k-1,1^{k+1})$ with degree $(+1,0,\ldots,0,-1),$ as shown below: \\
						{
\centering
\begin{tikzcd}[row sep=0.5cm, column sep=0.5cm]
							v_{0} \arrow[r, dashed] & v_{k} & v_{k+1} \arrow[l] & v_{n-k-2} \arrow[r] & v_{n-k-1} \arrow[r, dashed] & v_{n-1}.
						\end{tikzcd}
\captionof{figure}{}
}\item Subcase B.2 If $e_{n-k-1}$ is from $v_{n-k-1}$ to $v_{n-k-2}$, then $\text{td}p_{0,n-k-2}=-1$, as can be seen below: \\
						{
\centering
\begin{tikzcd}[row sep=0.5cm, column sep=0.5cm]
							v_{0} \arrow[r, dashed] & v_{k} & v_{k+1} \arrow[l] & v_{n-k-2} & v_{n-k-1} \arrow[r, dashed] \arrow[l] & v_{n-1}.
						\end{tikzcd}
\captionof{figure}{}
}\end{enumerate} The composition $(n-k-1,1^{k+1})$ in $\textbf{X}_{D} $ must correspond to a component of $n-k-1$ vertices with positive total degree. Accordingly, under case B, $c(n-k-1,1^{k+1})=2$ if $e_{n-k-1}$ is from $v_{n-k-2}$ to $v_{n-k-1}$; $c(n-k-1,1^{k+1})=1$ if $e_{n-k-1}$ is from $v_{n-k-1}$ to $v_{n-k-2}$.
		\end{quote}
		\end{enumerate}

		$\qedsymbol$
	\end{Proof}
	We can use the above proposition to determine some special oriented paths that are quite similar to directed paths.
	\begin{Corollary}\label{almost directed path}
		(1)	If $D$ is an oriented path with $2k+1$ vertices where $p_{0,k-1}$ is a directed path from $v_{0}$ to $v_{k-1}$, and $p_{k+1,2k}$ is a directed path from $v_{k+1}$ to $v_{2k}$, then $D$ can  be uniquely determined by $\textbf{X}_{D}$.\\
		(2)	If $D$ is an oriented path with $2k$ vertices where $p_{0,k-1}$ is a directed path from $v_{0}$ to $v_{k-1}$, and $p_{k,2k-1}$ is a directed path from $v_{k}$ to $v_{2k-1}$, then $D$ can  be uniquely determined by $\textbf{X}_{D}$.\\ 
	\end{Corollary}
	
        Proposition~\ref{right right} motivates us to consider how to determine the directions of edges deeper inside. We now change the direction of $e_{n-k}$ in the condition of Proposition~\ref{right right}.
	\begin{Proposition}\label{right left}
		Let $D$ be an oriented path with $n$ vertices where $p_{0,k}$ is a directed path from $v_{0}$ to $v_{k}$, $p_{n-k,n-1}$ is a directed path from $v_{n-k}$ to $v_{n-1}$, $e_{n-k}$ is from $v_{n-k}$ to $v_{n-k-1}$ and where $n> 2k+1$ with $k>1$. Then we can determine the direction of $e_{k+1}$ and $e_{n-k-1}$ from $\textbf{X}_{D}$.
	\end{Proposition}	\begin{Proof}
		The condition in this proposition is shown below: \\
		{
\centering
\begin{tikzcd}
			v_{0} \arrow[r, dashed] & v_{k}  & v_{n-k-1} & v_{n-k} \arrow[r, dashed] \arrow[l] & v_{n-1}.
		\end{tikzcd}
\captionof{figure}{}
}
		(The condition about $n$ and $k$ is to ensure that $n-k-1>k>1$)

		Consider the compositions in $\textbf{X}_{D}$ consisting of $k+1$ ``1's'' and one ``$n-k-1$''. The component of $n-k-1$ vertices must be $p_{i,i+n-k-2}$ for some $0\leq i\leq k+1$.
                \begin{enumerate}
                \item When $i=1$,  we find that  $\text{td} p_{1,n-k-1}=-2$, and that the corresponding composition is $(1^{k+1},n-k-1)$ with degree $(+2,+1,0,\ldots,0,-1,-2)$.
                \item In the case that  $2\leq i\leq k$, we determine that $\text{td} p_{i,i+n-k-2}=0$, and that the corresponding composition is $(1,n-k-1,1^{k})$ with degree $(+1,0,\ldots,0,-1)$.
                \item When $i=k+1$ and $e_{k+1}$ is from $v_{k}$ to $v_{k+1}$, 	the oriented path is shown below: \\
			{
\centering
\begin{tikzcd}
				v_{0} \arrow[r, dashed] & v_{k} \arrow[r] & v_{k+1}   & v_{n-k-1} & v_{n-k} \arrow[r, dashed] \arrow[l] & v_{n-1}.
			\end{tikzcd}
\captionof{figure}{}
}
 In this case, we determine that $\text{td} p_{k+1,n-1}=-1$, and that the corresponding composition is $(1^{k+1},n-k-1)$ with degree $(+1,0,\ldots,0,-1)$.\item When $i=k+1$ and $e_{k+1}$ is from $v_{k+1}$ to $v_{k}$, the oriented path is displayed below: \\
			{
\centering
\begin{tikzcd}
				v_{0} \arrow[r, dashed] & v_{k} & v_{k+1} \arrow[l]  & v_{n-k-1} & v_{n-k} \arrow[r, dashed] \arrow[l] & v_{n-1}
			\end{tikzcd}
\captionof{figure}{}
} In this case, we have that $\text{td} p_{k+1,n-1}=+1$, and that the corresponding composition is $(n-k-1,1^{k+1})$ with degree $(+1,+1,0,\ldots,0,-2)$.\item When $i=0$ and $e_{n-k-1}$ is from $v_{n-k-2}$ to $v_{n-k-1}$,	the oriented path is shown below: \\
			{
\centering
\begin{tikzcd}
				v_{0} \arrow[r, dashed] & v_{k}  & v_{n-k-2} \arrow[r] & v_{n-k-1} & v_{n-k} \arrow[r, dashed] \arrow[l] & v_{n-1}.
			\end{tikzcd}
\captionof{figure}{}
} In this case, we find that $\text{td} p_{0,n-k-2}=+1$, and that the corresponding composition is $(1,n-k-1,1^{k})$ with degree $(+2,+1,0,\ldots,0,-1,-2)$.\item When $i=0$ and $e_{n-k-1}$ is from $v_{n-k-1}$ to $v_{n-k-2}$, the oriented path is illustrated below: \\
			{
\centering
\begin{tikzcd}
				v_{0} \arrow[r, dashed] & v_{k} & v_{n-k-2} & v_{n-k-1} \arrow[l] & v_{n-k} \arrow[r, dashed] \arrow[l] & v_{n-1}.
			\end{tikzcd}
\captionof{figure}{}
} In this case, we can determine that $\text{td} p_{0,n-k-2}=-1$ and that the corresponding composition is $(1^{k},n-k-1,1)$ with degree $(+2,0,\ldots,0,-1,-1)$.\end{enumerate}
		In summary, $e_{k+1}$ is from $v_{k+1}$ to $v_{k}$ if and only if $c(n-k-1,1^{k+1})=1$; $e_{n-k-1}$ is from $v_{n-k-1}$ to $v_{n-k-2}$ if and only if $c(1^{k},n-k-1,1)=1$.
		$\qedsymbol$
	\end{Proof}
	For the cases obtained by reversing the direction of $e_{k}$ in the condition of Proposition~\ref{right right}, we consider not only the compositions consisting of one ``$n-k-1$'' and $k+1$ ``1's'' but also the compositions consisting of one ``$n-k-2$'' and $k+2$ ``1's''.
	\begin{Proposition}\label{left right}
		Let $D$ be an oriented path with $n$ vertices where $p_{0,k-1}$ is a directed path from $v_{0}$ to $v_{k-1}$, $p_{n-k-1,n-1}$ is a directed path from $v_{n-k-1}$ to $v_{n-1}$, $e_{k}$ is from $v_{k}$ to $v_{k-1}$, and $n> 2k+1$ with $k>1$. Then we can determine the direction of $e_{n-k-1}$ and $e_{k+1}$ from $\textbf{X}_{D}$.
	\end{Proposition}
	\begin{Proof}
		The condition in this proposition is shown below: \\
		{
\centering
\begin{tikzcd}
			v_{0} \arrow[r, dashed] & v_{k-1} & v_{k} \arrow[l] & v_{n-k-1} \arrow[r, dashed] & v_{n-1}.
		\end{tikzcd}
\captionof{figure}{}
}
		(The condition about $n$ and $k$ is to ensure that $n-k-1>n-k-2\geq k>1$)\\
		
		The first part of the proof is to discuss the direction of $e_{n-k-1}$. Consider whether the composition $(1^{k+2},n-k-2)$ appears in $\textbf{X}_{D}$. The component with $n-k-2$ vertices must be $p_{i,i+n-k-3}$ for some $0\leq i\leq k+2$, and its total degree must be negative.\begin{enumerate}
			\item When $i=0$, $\text{td} p_{0,n-k-3}\in \left\lbrace -1,+1\right\rbrace $. Even if $\text{td} p_{0,n-k-3}=-1$, ``$n-k-2$'' would not be the last part of the composition since $\text{td} v_{n-1}=-1$. \item When $2\leq i\leq k+1$,  we find that $\text{td} p_{i,i+n-k-3}=0$ for $2\leq i\leq k-1$ (If $k=2$, there are no such $i$.), $\text{td} p_{k,n-3}=+2$ and $\text{td} p_{k+1,n-2}\in \left\lbrace 0,+2\right\rbrace$ --- none of these would yield $(1^{k+2},n-k-2)$.\item When $i=k+2$, $\text{td} p_{k+2,n-1}\in \left\lbrace -1,+1\right\rbrace $. Even if $\text{td} p_{k+2,n-1}=-1$, ``$n-k-2$'' would not be the last part in the composition since $\text{td} v_{k-1}=-2$.

			\item	When $e_{n-k-1}$ is from $v_{n-k-2}$ to $v_{n-k-1}$, $\text{td} p_{1,n-k-2}=0$, the digraph is displayed below: \\
			{
\centering
\begin{tikzcd}
				v_{0} \arrow[r, dashed] & v_{k-1} & v_{k} \arrow[l] & v_{n-k-2} \arrow[r, dashed] & v_{n-1}.
			\end{tikzcd}
\captionof{figure}{}
}

			\item	When $e_{n-k-1}$ is from $v_{n-k-1}$ to $v_{n-k-2}$, $\text{td} p_{1,n-k-2}=-2$, the digraph is shown below: \\
			{
\centering
\begin{tikzcd}
				v_{0} \arrow[r, dashed] & v_{k-1} & v_{k} \arrow[l] & v_{n-k-2} & v_{n-k-1} \arrow[l] \arrow[r, dashed] & v_{n-1}.
			\end{tikzcd}
\captionof{figure}{}
} In this case, the corresponding composition is $(1^{k+2},n-k-2)$ with degree ($+1,+1,0\cdots,0,-2$).\end{enumerate}
		Accordingly, under the condition of this proposition, $(1^{k+2},n-k-2)$ appears in $\textbf{X}_{D}$ $\iff$ $e_{n-k-1}$ goes from $v_{n-k-1}$ to $v_{n-k-2}$.\\

Now we discuss $e_{k+1}$. The following proof is inspired by the proof of Proposition~\ref{right left}.	Consider the compositions in $\textbf{X}_{D}$ consisting of $k+1$ ``1's'' and one ``$n-k-1$''. The component of $n-k-1$ vertices must be $p_{i,i+n-k-2}$ for some $0\leq i\leq k+1$.
                \begin{enumerate}
                \item When $1\leq i\leq k-1$, we determine that $\text{td} p_{i,i+n-k-2}=0$, and the corresponding composition is $(1,n-k-1,1^{k})$ with degree $(+1,0,\ldots,0,-1)$.
                \item When $i=k$, we can determine $\text{td} p_{k,n-2}=+2$, and the corresponding composition is $(n-k-1,1^{k+1})$ with degree $(+2,+1,0,\ldots,0,-1,-2)$.
                \item 
			When $i=k+1$ and $e_{k+1}$ is from $v_{k}$ to $v_{k+1}$, the digraph is illustrated below: \\
			{
\centering
\begin{tikzcd}
				v_{0} \arrow[r, dashed] & v_{k-1} & v_{k} \arrow[l] \arrow[r] & v_{k+1}  & v_{n-k-1} \arrow[r, dashed] & v_{n-1}.
			\end{tikzcd}
\captionof{figure}{}
} In this case, we can determine that $\text{td} p_{k+1,n-1}=-1$ and that the corresponding composition is $(1^{k},n-k-1,1)$ with degree $(+2,+1,0,\ldots,0,-1,-2)$.\item When $i=k+1$ and $e_{k+1}$ is from $v_{k+1}$ to $v_{k}$, the digraph is shown below: \\
			{
\centering
\begin{tikzcd}
				v_{0} \arrow[r, dashed] & v_{k-1} & v_{k} \arrow[l] & v_{k+1}  \arrow[l] & v_{n-k-1} \arrow[r, dashed] & v_{n-1}.
			\end{tikzcd}
\captionof{figure}{}
}In this case, we can determine $\text{td} p_{k+1,n-1}=+1$, and that the corresponding composition is $(n-k-1,1^{k+1})$ with degree $(+1,+1,0,\ldots,0,-2)$.\item When $i=0$ and $e_{n-k-1}$ is from $v_{n-k-2}$ to $v_{n-k-1}$,	the digraph can be seen below: \\
			{
\centering
\begin{tikzcd}
				v_{0} \arrow[r, dashed] & v_{k-1} & v_{k} \arrow[l] & v_{n-k-2}  \arrow[r] & v_{n-k-1} \arrow[r, dashed] & v_{n-1}.
			\end{tikzcd}
\captionof{figure}{}
} In this case, we determine that $\text{td} p_{0,n-k-2}=+1$ and that the corresponding composition is $(n-k-1,1^{k+1})$ with degree $(+1,0,\ldots,0,-1)$.\item When $i=0$ and $e_{n-k-1}$ is from $v_{n-k-1}$ to $v_{n-k-2}$, the digraph is shown below:\\
			{
\centering
\begin{tikzcd}
				v_{0} \arrow[r, dashed] & v_{k-1} & v_{k} \arrow[l] & v_{n-k-2}  & v_{n-k-1} \arrow[r, dashed] \arrow[l] & v_{n-1}.
			\end{tikzcd}
\captionof{figure}{}
} In this case, we determine that $\text{td} p_{0,n-k-2}=-1$, and that the corresponding composition is $(1^{k},n-k-1,1)$ with degree $(+2,0,\ldots,0,-1,-1)$. \end{enumerate}
		Since we have determined the direction of $e_{n-k-1}$, we can determine $e_{k+1}$ according to $c(n-k-1,1^{k+1})$ in $\textbf{X}_{D}$. To be more precise, see the following:
		$$\begin{tabular}{c|c|c}
		\hline
		   $c(n-k-1)$ & $e_{n-k-1}$ is from $v_{n-k-2}$ & $e_{n-k-1}$ is from $v_{n-k-1}$\\
		    \hline
		   $e_{k+1}$ is from $v_k$       & 2 & 1 \\
		    \hline
		   $e_{k+1}$ is from $v_{k+1}$   & 3 & 2 \\\hline
		\end{tabular}$$
		$\qedsymbol$
	\end{Proof}The case when both $e_{k}$ and $e_{n-k}$ are changed is as follows: 
	\begin{Proposition}\label{left left}
		Let $D$ be an oriented path with $n$ vertices. Suppose that $p_{0,k-1}$ is a directed path from $v_{0}$ to $v_{k-1}$, $p_{n-k,n-1}$ is a directed path from $v_{n-k}$ to $v_{n-1}$, $e_{n-k}$ is directed from $v_{n-k}$ to $v_{n-k-1}$, $e_{k}$ is directed from $v_{k}$ to $v_{k-1}$, and $n>2k+1$ with $k>1$. Then we can determine the directions of $e_{k+1}$ and $e_{n-k-1}$ from $\textbf{X}_{D}$.
	\end{Proposition}
	\begin{Proof}
		The condition can be illustrated as follows: \\
		{
\centering
\begin{tikzcd}
			v_{0} \arrow[r, dashed] & v_{k-1} & v_{k} \arrow[l] & v_{n-k-1} & v_{n-k} \arrow[l] \arrow[r, dashed] & v_{n-1}.
		\end{tikzcd}
\captionof{figure}{}
}
(The condition about $n$ and $k$ ensures that $n-k-1>k>1$.)
Consider the compositions in $\textbf{X}_{D}$ consisting of $k+1$ ``1's'' and one ``$n-k-1$''. The component with $n-k-1$ vertices must be $p_{i,i+n-k-2}$ for some $0\leq i\leq k+1$.\begin{enumerate}
	
	 \item When $i=1$, we determine that $\text{td} p_{1,n-k-1}=-2$, and that the corresponding composition is $(1^{k+1},n-k-1)$ with degree $(+2,+1,0,\ldots,0,-1,-2)$.\item When $2\leq i\leq k-1$, we determine that $\text{td} p_{i,i+n-k-2}=0$, and that the corresponding composition is $(1,n-k-1,1^{k})$ with degree $(+1,0,\ldots,0,-1)$. If $k=2$, there are no such $i$. Hence, there are $k-2$ such $i$ whenever $ k\geq 2$.  \item When $i=k$, we determine that $\text{td} p_{k,n-2}=+2$, and that the corresponding composition is $(n-k-1,1^{k+1})$ with degree $(+2,+1,0,\ldots,0,-1,-2)$.\item When $i=k+1$ and $e_{k+1}$ is from $v_{k}$ to $v_{k+1}$, 	the digraph can be presented as follows: \\
		{
\centering
\begin{tikzcd}[row sep=0.5cm, column sep=0.4cm]
			v_{0} \arrow[r, dashed] & v_{k-1} & v_{k} \arrow[l] \arrow[r] & v_{k+1} & v_{n-k-1} & v_{n-k} \arrow[r, dashed] \arrow[l] & v_{n-1}.
		\end{tikzcd}
\captionof{figure}{}
} In this case, we can determine that $\text{td} p_{k+1,n-1}=-1$ and that the corresponding composition is $(1^{k},n-k-1,1)$ with degree $(+2,+1,0,\ldots,0,-1,-2)$.\item When $i=k+1$ and $e_{k+1}$ is from $v_{k+1}$ to $v_{k}$, the digraph can be displayed as follows:\\
		{
\centering
\begin{tikzcd}[row sep=0.5cm, column sep=0.4cm]
			v_{0} \arrow[r, dashed] & v_{k-1} & v_{k} \arrow[l] & v_{k+1} \arrow[l] & v_{n-k-1} & v_{n-k} \arrow[r, dashed] \arrow[l] & v_{n-1}.
		\end{tikzcd}
\captionof{figure}{}
} In this case, we deduce that $\text{td} p_{k+1,n-1}=+1$, and that the corresponding composition is $(n-k-1,1^{k+1})$ with degree $(+1,+1,0,\ldots,0,-2)$.\item When $i=0$ and $e_{n-k-1}$ is from $v_{n-k-2}$ to $v_{n-k-1}$,	the digraph can be observed as follows: \\
		{
\centering
\begin{tikzcd}[row sep=0.5cm, column sep=0.4cm]
			v_{0} \arrow[r, dashed] & v_{k-1} & v_{k} \arrow[l] & v_{n-k-2} \arrow[r] & v_{n-k-1} & v_{n-k} \arrow[r, dashed] \arrow[l] & v_{n-1}.
		\end{tikzcd}
\captionof{figure}{}
} In this case, we can determine $\text{td} p_{0,n-k-2}=+1$, and the corresponding composition is $(1,n-k-1,1^{k})$ with degree $(+2,+1,0,\ldots,0,-1,-2)$.\item When $i=0$ and $e_{n-k-1}$ is from $v_{n-k-1}$ to $v_{n-k-2}$, the digraph can be seen below: \\
		{
\centering
\begin{tikzcd}[row sep=0.5cm, column sep=0.4cm]
			v_{0} \arrow[r, dashed] & v_{k-1} & v_{k} \arrow[l] & v_{n-k-2} & v_{n-k-1} \arrow[l] & v_{n-k} \arrow[r, dashed] \arrow[l] & v_{n-1}.
		\end{tikzcd}
\captionof{figure}{}
} In this case, we can determine that $\text{td} p_{0,n-k-2}=-1$, and that the corresponding composition is $(1^{k},n-k-1,1)$ with degree $(+2,0,\ldots,0,-1,-1)$.
\end{enumerate}
		In summary, $e_{k+1}$ is from $v_{k+1}$ to $v_{k}$ if and only if $c(n-k-1,1^{k+1})=2$; $e_{n-k-1}$ is from $v_{n-k-2}$ to $v_{n-k-1}$ if and only if $c(1,n-k-1,1^{k})=k-1$.
		$\qedsymbol$
	\end{Proof}
	Combine Propositions~\ref{right right},~\ref{right left},~\ref{left right},~\ref{left left} and Corollary~\ref{almost directed path}, we have the following: 
	\begin{Corollary}\label{4 edges}
		Let $D$ be an oriented path with $n$ vertices where $p_{0,k-1}$ is a directed path from $v_{0}$ to $v_{k-1}$, $p_{n-k,n-1}$ is a directed path from $v_{n-k}$ to $v_{n-1}$, and $n-3>k>1$. Then we can determine the directions of $e_{k}$, $e_{k+1}$, $e_{n-k-1}$ and $e_{n-k}$ from $\textbf{X}_{D}$.
	\end{Corollary}
	\begin{Proof}
		We proceed by classification based on $n$.
		\begin{enumerate}
			\item
			When $n>2k+1$, we can first determine the directions of $e_{k}$ and $e_{n-k}$ using Proposition~\ref{right right} and then apply Proposition~\ref{right right},~\ref{right left},~\ref{left right} or~\ref{left left} depending on those directions to obtain the directions of $e_{k+1}$ and $e_{n-k-1}$.
                      \item	When $n=2k+1$ or $2k$, apply Corollary~\ref{almost directed path}.
                      \item
			When $n\leq 2k-1$, that is, $n-k\leq k-1$, the condition tells us the whole $D$.
                      \end{enumerate}
		
		Notice that	in our proof of Propositions~\ref{right left},~\ref{left right} and~\ref{left left}, we determine the desired directions of the edges using the compositions consisting of ``$n-k-1$'' and ``1's'' and those consisting of ``$n-k-2$'' and ``1's'' in $\textbf{X}_{D}$. Hence, we need that $n-k-2>1$, that is, $n-3>k.$
		$\qedsymbol$
	\end{Proof}
At the end of this subsection, we extend our results by replacing ``1'' with a general positive integer ``$i$''\footnote{The formation of this viewpoint benefits from offline discussions with Jindong Yan.}. For $n=mi$, an oriented path of $n$ vertices is called an \emph{i-directed path from $v_{0}$ to $v_{n-1}$} if $e_{ij}$ is from $v_{ij-1}$ to $v_{ij}$ for each $1\leq j\leq m-1$. Clearly, directed paths are 1-directed. Let $D$ be an $i$-directed path of $mi$ vertices. If we contract each $p_{ij,i(j+1)-1}$ of $D$ to one vertex $v'_{j}$ for $0\leq j\leq m-1$, then we turn $D$ into a directed path of $m$ vertices $D'$. Each edge $e_{j}$ in $D'$ corresponds to the edge $e_{ij}$ in $D$. The compositions where only one part is not $i$ in $\textbf{X}_{D}$ correspond to the compositions where only one part is not $1$ in $\textbf{X}_{D'}$. We have the following results:
\begin{Proposition}
If $D$ is an oriented path with $n=mi$ vertices, then $D$ is an $i$-directed path if and only if both $(i^{k},n-ki)$ and $(n-ki,i^{k})$ for each $k\leq \frac{m}{2}$ appear in $\textbf{X}_{D}$.
\end{Proposition}
		\begin{Proposition}\label{right right for i}
			Let $D$ be an oriented path with $n=mi$ vertices where $p_{0,ki-1}$ ($k\geq 1$) is an $i$-directed path from $v_{0}$ to $v_{ki-1}$, and $p_{n-ki,n-1}$ is an $i$-directed path from $v_{n-ki}$ to $v_{n-1}$ for some $1\leq k<m-2$. Then we can determine the direction of $e_{ki}$ and $e_{n-ki}$ from $\textbf{X}_{D}$.
		\end{Proposition}
		\begin{Proposition}\label{4 edge for i}
				Let $D$ be an oriented path with $n=mi$ vertices where $p_{0,ki-1}$ is an $i$-directed path from $v_{0}$ to $v_{ki-1}$, $p_{n-ki,n-1}$ is an $i$-directed path from $v_{n-ki}$ to $v_{n-1}$, and $m-3>k>1$. Then we can determine the directions of $e_{ki}$, $e_{(k+1)i}$, $e_{n-(k+1)i}$ and $e_{n-ki}$ from $\textbf{X}_{D}$.
			\end{Proposition}
			In fact, the conditions on $n$ in Propositions~\ref{right right for i} and~\ref{4 edge for i} can be weakened to $n>(k+2)i$ (to ensure that $n-(k+1)i>i$) and $n>(k+3)i$ (to ensure that $n-(k+1)i>i$ and $n-(k+2)i>i$) respectively, that is, $n$ does not have to be divisible by $i$. To see this, consider the compositions consisting of $k+1$ ``$i$'s'' and one ``$n-(k+1)i$'' and the compositions consisting of $k+2$ ``$i$'s'' and one ``$n-(k+2)i$''.
\subsection{Symmetry}\label{Symmetry}
In this section, we discuss the properties of CNSFs of oriented paths with respect to symmetry.
 For each composition $(a,b)$ of length 2 appearing in $\textbf{X}_{D}$, there are two possibilities --- either $e_{a}$ is from $v_{a-1}$	to $v_{a}$ or $e_{b}$ is from $v_{b}$	to $v_{b-1}$, which is illustrated as follows: \\
	{
\centering
\begin{tikzcd}
		v_{0} & \cdots & v_{a-1} \arrow[r] & v_{a}   & \cdots          & \cdots & v_{n} \\
		v_{0} & \cdots & \cdots            & v_{b-1} & v_{b} \arrow[l] & \cdots & v_{n}.
	\end{tikzcd}
\captionof{figure}{}
} When $D$ is ``symmetric'', the case is simple.
	

	\begin{Proposition}\label{symmetric}
		Let $D$ be an oriented path with $2k+1$ vertices. Suppose $D$ can be drawn as an axially symmetric shape with its axis of symmetry being the line passing through $v_{k}$ and perpendicular to $D$ (referred to as the \emph{central axis}). In that case, $D$ is uniquely determined by $\textbf{X}_{D}$.
	\end{Proposition}
	\begin{Proof}
		According to Proposition~\ref{degree and path sequence}, we can find that the underlying graph of $D$ is a path with $2k+1$ vertices from $\textbf{X}_{D}$. 
		 For $1\leq i\leq k$, we can find either $c(i,2k+1-i)=2$ or $c(2k+1-i,i)=2$ in $\textbf{X}_{D}$. If $c(i,2k+1-i)=2$, then $e_{i}$ is from $v_{j-1}$ to $v_{j}$ and $e_{2k-i+1}$ is from $v_{2k-i+1}$ to $v_{2k-i}$. If $c(2k+1-i,i)=2$, then $e_{i}$ is from $v_{j}$ to $v_{j-1}$ and $e_{2k-i+1}$ is from $v_{2k-i}$ to $v_{2k+1}$.
		$\qedsymbol$
	\end{Proof}
 The oriented paths satisfying the condition of Proposition~\ref{symmetric} are called \emph{symmetric}.	The condition of Proposition~\ref{symmetric} can be adjusted slightly.
	\begin{Proposition}\label{almost symmetric}
		Let $D$ be an oriented path with $2k+1$ vertices. Suppose $D$ can be changed into a symmetric oriented path by reversing only one edge. In that case, $D$ can be reconstructed from $\textbf{X}_{D}$.
	\end{Proposition}
	\begin{Proof}
		There is only one $i$ such that $c(i,2k+1-i)=c(2k+1-i,i)=1$ in $\textbf{X}_{D}$.
		For $1\leq j\leq k$ and $j\neq i$, either $c(j,2k+1-j)=2$ or $c(2k+1-j,j)=2$ in $\textbf{X}_{D}$, so we can determine the directions of $e_{j}$ and $e_{2k+1-j}$. For $e_{i}$ and $e_{2k+1-i}$, there are two cases: either both $e_{i}$ and $e_{2k+1-i}$ are from the vertex of the smaller index to the one of the larger index, or both $e_{i}$ and $e_{2k+1-i}$ are from the vertex of the larger index to the one of the smaller index --- the two cases of $D$ are isomorphic since we can reverse the labeling.
		$\qedsymbol$
	\end{Proof}
	We call an oriented path \emph{almost symmetric} if it satisfies the condition of either Proposition~\ref{symmetric} or Proposition~\ref{almost symmetric}. In other words, an oriented path with $2k+1$ vertices is almost symmetric if there are $k-1$ $j$'s in $\left[k\right] $ such that either $c(j,2k+1-j)=2$ or $c(2k+1-j,j)=2$ in $\textbf{X}_{D}$. Proposition~\ref{symmetric} and~\ref{almost symmetric} tell that almost symmetric oriented paths can be reconstructed from their CNSFs. There is also a similar result for oriented paths with an even number of vertices.
\begin{Proposition}\label{even symmetric}
	Let $D$ be an oriented path with $2k$ vertices. If $c(i,2k-i)=1$ holds only when $i=k$, then $D $ can be reconstructed from $\textbf{X}_{D}$.
\end{Proposition}
\begin{Proof}
	The condition tells us either $c(i,2k-i)=0$ or $c(i,2k-i)=2$ for  each $i\neq k$. Hence, we can determine the directions of all edges of $D$ except $e_{k}$. Notice that the two cases of $D$ where  $e_{k}$ is from either $v_{k-1}$ or $v_{k}$ are isomorphic.
	$\qedsymbol$
\end{Proof}
In this article, the \emph{distance} between two vertices in an oriented path is the distance between them in the underlying path.	
	\begin{Proposition}\label{part}
		Let $D$ be an oriented path whose leaves are both negative (positive). Then we can determine the number of sources (sinks) of $D$ from $\textbf{X}_{D}$. For each source (sink), we can also determine the distances between it and the leaves.
	\end{Proposition}
	\begin{Proof}
		According to Proposition~\ref{degree and path sequence}, we can find that the underlying graph of $D$ is a path of $n$ vertices. Every non-leaf vertex in $D$ is incident to 2 edges, so its total degree can be $+2$, 0, or $-2$.

 If $c(n-1,1)=2$  in $\textbf{X}_{D}$, then both leaves are negative. In this case,
		consider the compositions of length 3 whose first part is 1. The first part must correspond to a vertex whose total degree is positive, then it must be a non-leaf vertex whose total degree is $+2$, that is, it is a source. Notice that for every non-leaf source, deleting its incident edges causes one composition of length 3. So the number of sources is $\mathop{\sum}\limits_{a,b}c(1, a, b)$. For a source corresponding to $(1, a, b)$, the distances from it to the two leaves are $a$ and $b$, respectively.

		If $c(1,n-1)=2$ in $\textbf{X}_{D}$, then both leaves are positive. In this case, consider the compositions of length 3 whose last part is 1. The last part must correspond to a vertex whose total degree is negative, then it must be a non-leaf vertex whose total degree is $-2$, that is, it is a sink. Notice that for every non-leaf sink, deleting its incident edges causes one composition of length 3. So the number of sinks is $\mathop{\sum}\limits_{c,d}c(c,d,1)$. For a sink corresponding to  $(c,d,1)$, the distances from it to the two leaves are $c$ and $d$, respectively.
		$\qedsymbol$
	\end{Proof}
	Let $D$ be an oriented path of $n$ vertices. Suppose there are $s$ sources and $t$ sinks, then $|s-t|\leq 1$, and sinks and sources alternate. If we can determine the positions of all sources and sinks, then we can reconstruct $D$. Notice that positive leaves are sources, and negative ones are sinks.

	\begin{Corollary}\label{unimodal}
		(1)	Let $D$ be an oriented path whose leaves are both negative. If $D$ has exactly one source, then it can be uniquely determined by $\textbf{X}_{D}$.\\
		(2) Let $D$ be an oriented path whose leaves are both positive. If $D$ has exactly one sink, then it can be uniquely determined by $\textbf{X}_{D}$.
	\end{Corollary}\begin{Proof}
		(1)	According to Proposition~\ref{part}, there is only one composition of length 3 whose first part is 1 in $\textbf{X}_{D}$, and we let it be $(1, a, b)$. Let the vertices of $D$ be $v_{0},v_{1},\ldots,v_{a+b}$. Then the unique source is either $v_{a}$ or $v_{b}$. The direction of every edge is clear since it must belong to a directed path from the source to one leaf. Notice that the two cases --- where the source is $v_{a}$ and where it is$v_{b}$ --- yield isomorphic oriented paths, which are presented as follows: \\
		{
\centering
\begin{tikzcd}
			v_{0} &                                            & v_{a} \arrow[ll, dashed] \arrow[r, dashed] & v_{a+b}; \\
			v_{0} & v_{b} \arrow[l, dashed] \arrow[rr, dashed] &                                            & v_{a+b}.
		\end{tikzcd}
\captionof{figure}{}
}
		(2)	According to Proposition~\ref{part}, there is only one composition of length 3 whose last part is 1 in $\textbf{X}_{D}$, and we let it be $(c,d,1)$. Let the vertices of $D$ be $v_{0},v_{1},\ldots,v_{c+d}$. Then the unique sink is either $v_{c}$ or $v_{d}$. The direction of every edge is clear since it must belong to a directed path from one leaf to the sink. Notice that having the sink as $v_{c}$ or $v_{d}$ results in isomorphic oriented paths, as presented below: \\
		{
\centering
\begin{tikzcd}
			v_{0} \arrow[rr, dashed] &       & v_{c} & v_{c+d}; \arrow[l, dashed]  \\
			v_{0} \arrow[r, dashed]  & v_{d} &       & v_{c+d}. \arrow[ll, dashed]
		\end{tikzcd}
\captionof{figure}{}
}
		$\qedsymbol$
	\end{Proof}
An oriented path is called \emph{unimodal} if it satisfies the condition of either case (1) or case (2) in Corollary~\ref{unimodal}. Corollary~\ref{unimodal} shows that unimodal oriented paths can be reconstructed from their CNSFs. The result can be extended to the cases when there are odd sources or sinks.

Similar to the proof of Proposition~\ref{part}, we can obtain more information from the compositions of length 2 and 3 in $\textbf{X}_{D}$. For each composition $\alpha=(i,j)$ of length 2 appearing in $\textbf{X}_{D}$, there are two possibilities: either $\text{td} p_{0,i-1}$ or $\text{td} p_{j,n-1}$ is $+1$.

For $1\leq k< n$, if $c(n-k,k)=2$ in $\textbf{X}_{D}$, then $\text{td} p_{0,n-k-1}=\text{td} p_{k,n-1}=+1$, and $\text{td} p_{0,k-1}=\text{td} p_{n-k,n-1}=-1$. In this case, for each $(k,a,b)$ in $\textbf{X}_{D}$, since the total degree of the component of $k$ vertices can only be $+2$ --- it cannot be $+1$, we know that either $\text{td} p_{a,a+k-1}=+2$ or $\text{td} p_{b,b+k-1}=+2$. For each $(c,d,n-k)$ in $\textbf{X}_{D}$, since the total degree of the component of $n-k$ vertices can only be $-2$ --- it cannot be $-1$, we know that either $\text{td} p_{c,c+n-k-1}=-2$ or $\text{td} p_{d,d+n-k-1}=-2$.

For $1\leq i< n$, if there is no $(i,a,b)$ ($(a,b,i)$) appearing in $\textbf{X}_{D}$ for any $a,b$, then there is no $p_{j,j+i-1}$ whose total degree is $+2$ ($-2$).

For clarity, we will draw the sources above the sinks. If we contract $p_{i,j}$ to a vertex $v$, then $\text{td} p_{i,j}=+2$ $(-2)$ in $D$ if and only if $v$ is a convex (concave) point of the contracted digraph.

When leaves are both negative (positive) and the $k$ sources (sinks) are very ``close'' to each other, we can have a good result.
\begin{Proposition}\label{inflection-dense}
		Let $D$ be an oriented path with $n$ vertices $v_{0},\ldots,v_{n-1}$ and whose leaves are both negative (positive). If $D$ has exactly $k$ sources (sinks) $v_{i},v_{i+2},\ldots,v_{i+2(k-1)}$ where $2\leq i,i+2(k-1) \leq n-3$, then it can be uniquely determined by $\textbf{X}_{D}$. To be more precise,  $$\text{the negative leaf case holds}\iff \begin{cases}
		c(n-1,1)=2, \\
		c(n-2,2)=2, \\
		\mathop{\sum}\limits_{a,b}c(1, a, b)=k,\\
		\mathop{\sum}\limits_{s,r}c(2,s,r)=2;
		\end{cases}$$ while $$\text{the positive leaf case holds} \iff \begin{cases}c(1,n-1)=2,\\c(2,n-2)=2,\\\mathop{\sum}\limits_{a,b}c(a, b,1)=k,\\\mathop{\sum}\limits_{s,r}c(s,r,2)=2;\\\end{cases} 
		$$and when one of the cases holds, $D$ can be determined by  $\min\left\lbrace b|\exists a, c(1,a,b)\neq0\right\rbrace$ ($\min\left\lbrace b|\exists a, c(a,b,1)\neq0\right\rbrace$) and $k$.
	\end{Proposition}
	\begin{Proof}
We first consider the case where both leaves are negative. For the positive leaf case, refer to the adjusted statements marked with $\dagger$ below.
According to Proposition~\ref{part}, there are exactly $k$ sources $\dagger$ in $D$ $ \iff$ $\mathop{\sum}\limits_{a,b}c(1, a, b)\dagger=k.$ Let the indices of the $k$ sources $\dagger$ be $j_{1}<j_{2}<\cdots<j_{k}$. Then there is a sink $\dagger$ $v_{s_{t}}$ between each $v_{j_{t}}$ and $v_{j_{t+1}}$, as can be seen below:\\
		{
\centering
\begin{tikzcd}
			& v_{j_{1}} \arrow[ld, dashed] \arrow[rd, dashed] &                                                 & v_{j_{2}} \arrow[ld, dashed] \arrow[rd, dashed] &                                              & v_{j_{k}} \arrow[ld, dashed] \arrow[rd, dashed] &                              \\
			v_{0}                     &                                                 & v_{s_{1}}                                       &                                                 & \cdots                                       &                                                 & v_{n-1}                      \\
			(v_{0} \arrow[rd, dashed] &                                                 & v_{s_{1}} \arrow[ld, dashed] \arrow[rd, dashed] &                                                 & \cdots \arrow[ld, dashed] \arrow[rd, dashed] &                                                 & v_{n-1}\dagger). \arrow[ld, dashed] \\
			& v_{j_{1}}                                       &                                                 & v_{j_{2}}                                       &                                              & v_{j_{k}}                                       &                             
		\end{tikzcd}
\captionof{figure}{}
}

	 For any $p_{a,b}$ whose total degree is +2 $\dagger$, the vertex $v_{a}$ lies between a sink $\dagger$ $v_{i}$ where $i<a$ and a source $\dagger$ $v_{j}$ where $a\leq j$, while $v_{b}$ lies between a source $\dagger$ $v_{s}$ where $s\leq b$ and a sink $\dagger$ $v_{t}$ where $b<t$. This can be seen as follows:\\
		{
\centering
\begin{tikzcd}
			& v_{j} \arrow[ld, "v_{a}", dashed] & v_{s} \arrow[rd, "v_{b}", dashed] &                                     \\
			v_{i}                              &                                   &                                   & v_{t}                               \\
			(v_{i} \arrow[rd, "v_{a}", dashed] &                                   &                                   & v_{t} \dagger ). \arrow[ld, "v_{b}", dashed] \\
			& v_{j}                             & v_{s}                             &                                    
		\end{tikzcd}
\captionof{figure}{}
}
		When $c(n-2,2)=2$ $\dagger$, there is no component with 2 vertices of total degree $+1$ $\dagger$; thus, for the compositions of the form $(2,c,d)$ $\dagger$ appearing in $\textbf{X}_{D}$, the total degree of the first $\dagger$ part must be $+2$ $\dagger$. If $\text{td}p_{a,a+1}=+2$ $\dagger$, then $a$ can be $j_{1}-1, j_{1}, j_{2}-1, j_{2}, \ldots, j_{k}-1, j_{k}$. For each $t$, $\text{td}p_{j_{t},j_{t}+1}=+2$ $\dagger$ $\iff$ $s_{t}>j_{t}+1$; $\text{td}p_{j_{t}-1,j_{t}}=+2$ $\dagger$ $\iff$ $j_{t}>s_{t-1}+1$. Therefore, $\mathop{\sum}\limits_{s,r}c(2,s,r)=2$ $\dagger$ in $\textbf{X}_{D}$ $\iff$ $s_{t}=j_{t}+1$ and $j_{t+1}=s_{t}+1$, that is, $j_{t+1}=j_{t}+2$, for all $1\leq t\leq k-1$.

Suppose that $j_{t+1}=j_{t}+2$ for each $1\leq t\leq k-1$ from $\textbf{X}_{D}$. Let the $2k$ indices from which $j_{1},j_{2},\ldots,j_{k}$ can take values be $i_{1}\leq i_{2}\leq\cdots\leq i_{k}\leq n-1-i_{k}\leq\cdots \leq n-1-i_{1}$ --- these values are obtained by Proposition~\ref{part}, and in particular, $i_{1}=\min\left\lbrace b|\exists a, c(1,a,b)\dagger\neq0\right\rbrace.$ Since either $i_{1}$ or $n-1-i_{1}$ must be chosen as the index of the source $\dagger$ closest to a leaf, the sequence ($j_{1}<j_{2}<\cdots<j_{k}$) can only be either $(i_{1}<i_{1}+2<\cdots<i_{1}+2(k-1))$ or $(n-3-i_{1}-2k<n-1-i_{1}-2k<\cdots<n-1-i_{1})$ --- the two cases are isomorphic.

For the positive leaf case, replace all $\dagger$-marked elements as follows:
		\begin{itemize}
		\item ``$c(1, a, b)$'' $\to$ ``$c(a,b,1)$''
		\item ``$c(n-2,2)=2$'' $\to$ ``$c(2,n-2)=2$''
			\item ``first'' $\to$ ``last''
			\item ``sources'' $\to$ ``sinks''
			\item ``($2,c,d$)'' $\to$ ``($c,d,2$)''
			\item ``$+1$'' $\to$ ``$-1$''
			\item ``$+2$'' $\to$ ``$-2$''
			\item ``source'' $\to$ ``sink''
			\item ``sink'' $\to$ ``source''
			\item ``$\mathop{\sum}\limits_{c,d}c(2,s,r)=2$'' $\to$ ``$\mathop{\sum}\limits_{c,d} c(s,r,2)=2$''
		\end{itemize}
		For illustrations, see as in ``()''.
		$\qedsymbol$
	\end{Proof}
	Oriented paths satisfying the condition of either case (1) or case (2) in Proposition~\ref{inflection-dense} are called \emph{inflection-dense}.
\subsection{$n$ vertices}
In this section, we prove that oriented paths with up to 24 vertices can be determined by their CNSFs. For $1\leq n\leq 7$, we prove this using theoretical methods; for $8\leq n\leq 24$, we verify it through program calculations.
	\begin{Proposition}
		Let $D$ be an oriented path with $n$ vertices. If $n\leq 5$, then $D$ can be uniquely determined by $\textbf{X}_{D}$.
	\end{Proposition}
	\begin{Proof}
		According to Proposition~\ref{degree and path sequence}, we can determine that the underlying graph of $D$ is a path with $n$ vertices from $\textbf{X}_{D}$.\begin{enumerate}
			\item	When $n=1$, $D$ is the trivial graph.
			\item	When $n=2$, $D$ is a directed path of length 1.
			\item	When $n=3$, $D$ is almost symmetric.
			\item	When $n=4$, the underlying graph of $D$ can be seen as a double star $S_{1,1}$; this case has been discussed as Case D in the proof of Proposition~\ref{double star}.
	
	\item When $n=5$, we discuss the classification based on the condition of the leaves.

			\begin{itemize}
				\item	If one leaf of $D$ is positive and the other one is negative, denote the vertices from the positive leaf to the negative one by $v_{0},\ldots,v_{4}$. Then $D$ satisfies the condition of case (1) in Proposition~\ref{right right}.\item	If both leaves are positive or negative, then $D$ is almost symmetric.\end{itemize}\end{enumerate}
			$\qedsymbol$
		\end{Proof}
		The cases when $n=6$ are much more complicated.
		\begin{Proposition}\label{n=6}
			Let $D$ be an oriented path with 6 vertices. Then it can be uniquely determined by $\textbf{X}_{D}$.
		\end{Proposition}
		\begin{Proof}
			We discuss the classification based on the condition of the leaves. 
			\begin{enumerate}
				\item	We first consider the case that both leaves are positive or negative. \begin{quote} When denoting the vertices by $v_{0},\ldots,v_{5}$ and edges by $e_{1},\ldots,e_{5}$, we require that $e_{3}$ is from $v_{2}$ to $v_{3}$ --- this constraint on the labeling helps avoid redundant discussions about isomorphism classes.
					\begin{enumerate}\item	When both leaves are positive, which happens if and only if $c(1,5)=2$ in $\textbf{X}_{D}$, we only need to know the directions of $e_{2}$ and $e_{4}$. The 4 possibilities are as follows: \\
						{
\centering
\begin{tikzcd}
							(1)v_{0} \arrow[r] & v_{1} \arrow[r] & v_{2} \arrow[r]           & v_{3} \arrow[r] & v_{4}           & v_{5}; \arrow[l] \\
							(2)v_{0} \arrow[r] & v_{1} \arrow[r] & v_{2} \arrow[r]           & v_{3}           & v_{4} \arrow[l] & v_{5}; \arrow[l] \\
							(3)v_{0} \arrow[r] & v_{1}           & v_{2} \arrow[r] \arrow[l] & v_{3} \arrow[r] & v_{4}           & v_{5}; \arrow[l] \\
							(4)v_{0} \arrow[r] & v_{1}           & v_{2} \arrow[r] \arrow[l] & v_{3}           & v_{4} \arrow[l] & v_{5}; \arrow[l]
						\end{tikzcd}
\captionof{figure}{}
}
					
						\begin{itemize}
							\item Cases (1) and (2) are unimodal. \item Case (3) satisfies the condition of Proposition~\ref{even symmetric}.
						\item 	If $\textbf{X}_{D}$ satisfies $c(1,5)=2$, but it is not case (1), (2), or (3), then it is case (4). \end{itemize} \item	When both leaves are negative, which happens if and only if $c(5,1)=2$ $\textbf{X}_{D}$, we only need to know the directions of $e_{2}$ and $e_{4}$. The 4 possibilities are as follows: \\
						{
\centering
\begin{tikzcd}
							(5)v_{0} & v_{1} \arrow[l] \arrow[r] & v_{2} \arrow[r]           & v_{3} \arrow[r] & v_{4} \arrow[r]           & v_{5}; \\
							(6)v_{0} & v_{1} \arrow[l] \arrow[r] & v_{2} \arrow[r]           & v_{3}           & v_{4} \arrow[r] \arrow[l] & v_{5}; \\
							(7)v_{0} & v_{1} \arrow[l]           & v_{2} \arrow[r] \arrow[l] & v_{3} \arrow[r] & v_{4} \arrow[r]           & v_{5}; \\
							(8)v_{0} & v_{1} \arrow[l]           & v_{2} \arrow[r] \arrow[l] & v_{3}           & v_{4} \arrow[r] \arrow[l] & v_{5}.
						\end{tikzcd}
\captionof{figure}{}
}
						\begin{itemize}
							\item
							Cases (5) and (7) are unimodal.\item Case (6) satisfies the condition of Proposition~\ref{even symmetric}.\item If $\textbf{X}_{D}$ satisfies $c(5,1)=2$, but it is not case (5), (6), or (7), then it is case (8). \end{itemize}\end{enumerate}\end{quote}
				\item	When one leaf is positive and the other one is negative, denote the vertices from the positive leaf to the negative one by $v_{0},\ldots,v_{5}$ and edges by $e_{1},\ldots,e_{5}$ to avoid isomorphism. According to Corollary~\ref{4 edges}, we can determine the direction of $e_{2}$, $e_{3}$, and $e_{4}$ by $\textbf{X}_{D}$. \end{enumerate}
			$\qedsymbol$
		\end{Proof}
		\begin{Proposition}\label{n=7}		
		Let $D$ be an oriented path with 7 vertices, then it can be uniquely determined by $\textbf{X}_{D}$.
              \end{Proposition}
              \begin{Proof}
                We discuss by classification based on $c(6,1)$.
                \begin{enumerate}
                \item Case A:  $c(6,1)=0$ \begin{quote}
											In this case, both leaves are positive. Now we consider  $c(5,2)$.
											\begin{enumerate}
												\item Subcase A.1:  If $c(5,2)=0$ or 2, then 
												 $D$ is almost symmetric.
												
												\item Subcase A.2:  If $c(5,2)=1$, denote the edges by $v_{0},\ldots,v_{6}$ and the edges by $e_{1},\ldots,e_{6}$ such that $e_{2}$ is from $v_{1}$ to $v_{2}$.
												
												Further, under Subcase A.2: 
												\begin{itemize}
													\item Subsubcase A.2.1:  If $c(4,3)=0$ or 2, then $D$ is almost symmetric.
													\item Subsubcase A.2.2:  If $c(4,3)=1$, then there are two cases: 
													{
\centering
\begin{tikzcd}[row sep=0.5cm, column sep=0.4cm]
														A.2.2.1: v_{0} \arrow[r] & v_{1} \arrow[r] & v_{2} \arrow[r] & v_{3} \arrow[r] & v_{4} \arrow[r]           & v_{5} & v_{6}\text{;} \arrow[l] \\
														A.2.2.2: v_{0} \arrow[r] & v_{1} \arrow[r] & v_{2}           & v_{3} \arrow[l] & v_{4} \arrow[r] \arrow[l] & v_{5} & v_{6}\text{.} \arrow[l]
													\end{tikzcd}
\captionof{figure}{}
}
													Case A.2.2.1 is unimodal while case A.2.2.2 is not, so they can be distinguished from each other under the condition of Subsubcase A.2.2.
                                                                                                      \end{itemize}
                                                                                              \end{enumerate}
										\end{quote}
										\item Case B: $c(6,1)=1$ \begin{quote}
											In this case, one leaf is positive while the other one is negative. Let $k=2$ and denote the positive leaf by $v_{0}$. We can determine the directions of all edges of $D$ according to Corollary~\ref{4 edges}. 
										\end{quote}
										
										\item Case C: $c(6,1)=2$\begin{quote}
											In this case, both leaves are negative. Now we consider  $c(5,2)$.
											\begin{enumerate}
												\item Subcase C.1:  If $c(5,2)=0$ or 2, then $D$ is almost symmetric.
												\item Subcase C.2:  If $c(5,2)=1$, denote the edges by $v_{0},\ldots,v_{6}$ and the edges by $e_{1},\ldots,e_{6}$ such that $e_{2}$ is from $v_{1}$ to $v_{2}$.
												\begin{itemize}
													\item Subsubcase C.2.1:  If $c(4,3)=0$ or 2, then $D$ is almost symmetric.
													\item Subsubcase C.2.2:  If $c(4,3)=1$, then there are two cases: 
													
												{
\centering
\begin{tikzcd}[row sep=0.5cm, column sep=0.4cm]
																								C.2.2.1:v_{0} & v_{1}  \arrow[l] \arrow[r] & v_{2} \arrow[r] & v_{3} \arrow[r] & v_{4} \arrow[r] & v_{5} \arrow[r] & v_{6}; \\
																						C.2.2.2:v_{0} & v_{1} \arrow[l] \arrow[r]  & v_{2} \arrow[r] & v_{3} \arrow[r] & v_{4} \arrow[r] & v_{5} \arrow[r] & v_{6}.																					
																					\end{tikzcd}
\captionof{figure}{}
}
													Case C.2.2.1 is unimodal while case C.2.2.2 is not, so they can be distinguished from each other under the condition of Subsubcase C.2.2.
												\end{itemize}
											\end{enumerate}
										\end{quote}
									\end{enumerate}
									$\qedsymbol$
									
								\end{Proof}
		Proposition~\ref{n=7} suggests the following:
		\begin{Proposition}	
				Let $D$ be an oriented path with $n\geq 7$ vertices, then $e_{1},e_{2},e_{3},e_{n-3},e_{n-2}$, and $e_{n-1}$ (under our rule of labeling) can be reconstructed  by $\textbf{X}_{D}$.
			\end{Proposition}
			\begin{Proof}
							We discuss by classification based on $c(n-1,1)$.
							\begin{enumerate}
								\item Case A:  $c(n-1,1)=0$ \begin{quote}
									In this case, we know $e_{1}$ and $e_{n-1}$. Now we consider  $c(n-2,2)$.
									\begin{enumerate}
										\item Subcase A.1:  If $c(n-2,2)=0$ or 2, then $e_{2}$ and $e_{n-2}$ can be determined by $\textbf{X}_{D}$.

										 \begin{itemize}
											\item Subsubcase A.1.1:  If $c(n-3,3)=0$ or 2, then $e_{3}$ and $e_{n-3}$ can be determined by $\textbf{X}_{D}$.
											\item Subsubcase A.1.2:  If $c(n-3,3)=1$, then we denote the edges by $v_{0},\ldots,v_{n-1}$ and the edges by $e_{1},\ldots,e_{n-1}$ such that $e_{3}$ is from $v_{2}$ to $v_{3}$. In this case, we have that $e_{n-3}$ is from $v_{n-4}$ to $v_{n-3}$.
										\end{itemize}
										
										\item Subcase A.2:  If $c(n-2,2)=1$, denote the edges by $v_{0},\ldots,v_{n-1}$ and the edges by $e_{1},\ldots,e_{n-1}$ such that $e_{2}$ is from $v_{1}$ to $v_{2}$.
										
										Further, under Subcase A.2: 
										\begin{itemize}
											\item Subsubcase A.2.1:  If $c(n-3,3)=0$ or 2, then $e_{3}$ and $e_{n-3}$ can be determined by $\textbf{X}_{D}$.
											\item Subsubcase A.2.2:  If $c(n-3,3)=1$, then there are two cases: 
											\begin{footnotesize}{
\centering
\begin{tikzcd}[row sep=0.3cm, column sep=0.3cm]
												(1): v_{0} \arrow[r] & v_{1} \arrow[r] & v_{2} \arrow[r] & v_{3}\ldots v_{n-4} \arrow[r] & v_{n-3} \arrow[r]           & v_{n-2} & v_{n-1}; \arrow[l] \\
												(2): v_{0} \arrow[r] & v_{1} \arrow[r] & v_{2}           & v_{3}\ldots v_{n-4} \arrow[l] & v_{n-3} \arrow[r] \arrow[l] & v_{n-2} & v_{n-1}. \arrow[l]
											\end{tikzcd}
\captionof{figure}{}
}\end{footnotesize}
										In	case (1) $c(1^{3},n-3)=2$ --- the component of $n-3$ vertices can be $p_{2,n-2}$ or $p_{3,n-1}$, while in case (2), $c(1^{3},n-3)=1$ --- the component of $n-3$ vertices can only be $p_{2,n-2}$.
										\end{itemize}
									\end{enumerate}
								\end{quote}
								\item Case B: $c(n-1,1)=1$ \begin{quote}
									In this case, one leaf is positive while the other one is negative. Let $k=2$ and denote the positive leaf by $v_{0}$. We can determine the directions of $e_{2},e_{3},e_{n-3}$, and $e_{n-2}$ according to Corollary~\ref{4 edges}. 
								\end{quote}
								
								\item Case C: $c(n-1,1)=2$\begin{quote}
									In this case, both leaves are negative. Now we consider  $c(n-2,2)$.
									\begin{enumerate}
										\item Subcase C.1:  If $c(n-2,2)=0$ or 2, then $e_{2}$ and $e_{n-2}$ can be determined by $\textbf{X}_{D}$.
											\item Subsubcase C.1.1:  If $c(n-3,3)=0$ or 2, then $e_{3}$ and $e_{n-3}$ can be determined by $\textbf{X}_{D}$.
											\item Subsubcase C.1.2:  If $c(n-3,3)=1$, then we denote the edges by $v_{0},\ldots,v_{n-1}$ and the edges by $e_{1},\ldots,e_{n-1}$ such that $e_{3}$ is from $v_{2}$ to $v_{3}$. In this case, we have that $e_{n-3}$ is from $v_{n-4}$ to $v_{n-3}$.
										
										\item Subcase C.2:  If $c(n-2,2)=1$, denote the edges by $v_{0},\ldots,v_{n-1}$ and the edges by $e_{1},\ldots,e_{n-1}$ such that $e_{2}$ is from $v_{1}$ to $v_{2}$.
										\begin{itemize}
											\item Subsubcase C.2.1:  If $c(n-3,3)=0$ or 2, then $e_{3}$ and $e_{n-3}$ can be determined by $\textbf{X}_{D}$.
											\item Subsubcase C.2.2:  If $c(n-3,3)=1$, then there are two cases: 
											\begin{footnotesize}{
\centering
\begin{tikzcd}[row sep=0.3cm, column sep=0.3cm]
											(3):v_{0} & v_{1}  \arrow[l] \arrow[r] & v_{2} \arrow[r] & v_{3}\ldots v_{n-4} \arrow[r] & v_{n-3} \arrow[r] & v_{n-2} \arrow[r] & v_{n-1}; \\
											(4):v_{0} & v_{1} \arrow[l] \arrow[r]  & v_{2} \arrow[r] & v_{3}\ldots v_{n-4} \arrow[r] & v_{n-3} \arrow[r] & v_{n-2} \arrow[r] & v_{n-1}.
											\end{tikzcd}
\captionof{figure}{}
}\end{footnotesize}
											In	case (3) $c(n-3,1^{3})=2$ --- the component of $n-3$ vertices can be $p_{0,n-4}$ or $p_{1,n-3}$, while in case (4), $c(n-3,1^{3})=1$ --- the component of $n-3$ vertices can only be $p_{1,n-3}$.
										\end{itemize}
									\end{enumerate}
								\end{quote}
							\end{enumerate}
							$\qedsymbol$
							
						\end{Proof}
						
		When $n=2j$, there are $2^{2j-1}$ possibilities considering the sequences of edges, and each case $D$ is isomorphic to the one obtained by reversing the labeling of $D$. Hence, there are $2^{2j-2}$ non-isomorphic oriented paths of $2j$ vertices. 

When $n=2j+1$, there are $2^{2j}$ possibilities considering the sequences of edges, and each non-symmetric case $D$ is isomorphic to the one obtained by reversing the labeling of $D$. There are $2^{j}$ symmetric oriented paths. Hence, there are $$2^{2j}-(2^{2j}-2^{j})/2=2^{2j}-2^{2j-1}+2^{j-1}=2^{2j-1}+2^{j-1}$$ non-isomorphic oriented paths of $2j+1$ vertices.

For $n>7$, the number of non-isomorphic oriented paths might be too large to calculate one by one. When discussing the oriented paths with $n$ vertices, we often consider the compositions where one part is ``$k$'' and the other parts are all ``1's'' for some $1<k<n$. This motivates us to write an algorithm to calculate all compositions appearing in $\textbf{X}_{D}$ where there is exactly one part larger than 1 (except for the composition $(n)$). The main steps are as follows:
\begin{enumerate}
\item Step 1: Generate all non-isomorphic oriented paths of $n$ vertices. \begin{quote}
For each oriented path, we use a sequence of length $n-1$ consisting of ``1's'' and ``$-1$'s'' to represent the sequence of edges, where ``1'' means the edge is from left to right, and ``$-1$'' means the other direction.
\begin{enumerate}
\item When $n$ is even, there are an odd number of edges, and we require the middle edge to go from left to right to avoid isomorphism.
\item When $n$ is odd, there are an even number of edges. For a non-symmetric oriented path $D$, we only keep either $D$ or the oriented path obtained by reversing the labeling of $D$. Since we have already handled symmetric oriented paths, we skip them.
\end{enumerate}\end{quote}
\item Step 2: For each oriented path $D$ obtained in step 1, calculate all compositions appearing in $\textbf{X}_{D}$ where there is exactly one part larger than 1. We get a list $[m_{2},m_{3},\ldots,m_{n-1}]$ where $m_{k}$ is the multiset of compositions appearing in $\textbf{X}_{D}$ consisting of one ``$k$'' and $n-k$ ``1's''.
\item Step 3: Check whether the sequences obtained in step 2 distinguish non-isomorphic oriented paths.
\end{enumerate}
We offer the Python implementation of the program in Appendix~\ref{app sec A}.
For $3\leq n \leq 24$, running ``main(n)'', we get ``Success! This method can distinguish oriented paths of $n$ vertices.'' As $n$ increases, the computational complexity of the program grows exponentially. A 16 GB-memory computer can handle this program at a maximum of approximately $n=21$, while calculating $n=24$ may require a 256 GB-memory computer. When running ``main (25)'' on the 256 GB-memory computer we rented, it displayed the error message: ``Out of memory resources, unable to process this command.''
Readers with better devices can try larger values of $ n$.To handle the computational complexity of the algorithm, we implemented a driver script (\url{maxn.py}) that monitors memory usage during execution. The driver estimates memory consumption incrementally and terminates the process before potential overflow occurs, ensuring stable performance across different hardware configurations. The core mathematical functions are written in \url{cnsf_paths.py} (same code as in the Appendix~\ref{app sec A}), which is imported by the driver at runtime.

Now we have that CNSFs can distinguish non-isomorphic oriented paths with up to 24 vertices. We propose the following conjecture:
		\begin{Conjecture}
				Non-isomorphic oriented paths have different CNSFs.
			\end{Conjecture}
\appendix
\section{Python code}\label{app sec A}
			    \begin{lstlisting}[language=Python, caption={Python code}, label={python}]
		
from itertools import product
import time
#Generate all oriented paths of n vertices. A sequence of length n-1 consisting of "1"s and "-1"s represents the sequence of edges, where "1" means the edge is from the left to the right, and "-1" means the other direction.
def generate_sequences(n):
    seq_length = n - 1
    return product([1, -1], repeat=seq_length) if seq_length > 0 else []
# Calculate the isomorphic digraph of an oriented path.
def iso(path_a):
    return tuple(-path_a[len(path_a)-1-i] for i in range(len(path_a)))
#Check whether an oriented path is symmetric
def sym(path_a):
    return all(path_a[i] + path_a[len(path_a)-1-i] == 0 for i in range(len(path_a)//2))
#Fix an oriented path, and integers i and k, calculate the composition corresponding to p_{i,i+k-1}. "td_values" and "prefix_gt" are from "path".
def one_composition(path, td_values, prefix_gt, i, k):
    n_vertices = len(td_values)
    if i == 0:
        td_p_i = path[k - 1]
    elif i == n_vertices - k:
        td_p_i = -path[n_vertices - k - 1]
    else:
        td_p_i = path[i + k - 1] - path[i - 1]
    start, end = i, i + k - 1
    t = td_p_i
    ones_before_k = prefix_gt[t][start] + (prefix_gt[t][n_vertices] - prefix_gt[t][end + 1])
    ones_after_k = n_vertices - k - ones_before_k
    return tuple([1] * ones_before_k + [k] + [1] * ones_after_k)
def compositions_with_ones(path):
    n_vertices = len(path) + 1
    td_values = [0] * n_vertices #We use td_values to record the total degree of vertices.
    if n_vertices >= 1:
        td_values[0] = path[0]
    if n_vertices >= 2:
        td_values[-1] = -path[-1]
    for j in range(1, n_vertices - 1):
        td_values[j] = path[j] - path[j - 1]
  
    prefix_gt = {}
    for t in [-2, -1, 0, 1, 2]:
        prefix = [0] * (n_vertices + 1)
        for j in range(n_vertices):
            prefix[j + 1] = prefix[j] + (1 if td_values[j] > t else 0)
        prefix_gt[t] = prefix #prefix_gt[t][i] is the number of vertices whose total degrees are larger than t and whose indices are smaller than i. We use prefix_gt to reduce the number of comparisons.
    result = []
    for k in range(2, n_vertices):
    #Fix an oriented path and an integer k, get the list k_compositions whose i-th element is the composition corresponding to p_{i,i+k-1}.
        max_i = n_vertices - k + 1
        k_compositions = []
        for i in range(max_i):
            k_compositions.append(one_composition(path, td_values, prefix_gt, i, k))
        result.append(tuple(sorted(k_compositions)))
    return tuple(result)
    
#When n is odd, delete the isomorphic oriented paths and check whether there are non-isomorphic oriented paths with the same result of "compositions_with_ones". For a non-symmetric oriented path D, we keep either D or iso(D).
def process_sequences(sequences, check_symmetry):
    processed = set()#Collect the sequences we have 
    seen_results = set()
    for seq in sequences:
        seq_tuple = tuple(seq)
        #We have proved for symmetric oriented paths, so we skip them.
        if check_symmetry and sym(seq_tuple):
            continue
        if seq_tuple in processed:
            continue
        result = compositions_with_ones(seq_tuple)
        if result in seen_results:
            return False
        seen_results.add(result)
        processed.add(seq_tuple)
        processed.add(iso(seq_tuple))
    return True

def main(n):
    print(f"Begin for oriented paths of {n} vertices")
    if n <= 2:
        print(f"There is only one oriented path of {n} vertices")
        return True
    if n % 2 == 1:
        sequences = generate_sequences(n)
        if not process_sequences(sequences, check_symmetry=True):
            print(f"Failed. This method can not distinguish oriented paths of {n} vertices")
            return False
    else:#When n is even, there are an odd number of edges, and we require the middle edge to go from left to right to avoid isomorphism.
        i = (n-2) // 2        
        pre_sequences = list(product([1, -1], repeat=i))
        post_sequences = list(product([1, -1], repeat=i))                
        seen_results = set()         
        for pre in pre_sequences:
            for post in post_sequences:           
                seq_tuple = pre + (1,) + post                      
                result = compositions_with_ones(seq_tuple)
                if result in seen_results:
                    print(f"Failed. This method can not distinguish oriented paths of {n} vertices")
                    return False
                seen_results.add(result)
    print(f"Success! This method can distinguish oriented paths of {n} vertices.")
    return True

\end{lstlisting}
\section{Data Availability Statement}
The Python code developed for this study is available at Zenodo \cite{zenodo2026code}.The Python code developed for this study is available at Zenodo via the following DOI: https://doi.org/10.5281/zenodo.22652747. The repository contains two files: \url{cnsf_paths.py}, which implements the core mathematical functions, and \url{maxn.py}, a driver script for testing and demonstration.

	%

\begin{thebibliography}{10}
		\bibitem{zenodo2026code}
		Lingxiao Hao, and Shenglin Zhu. (2026) \emph{Code for [The chromatic noncommutative symmetric function of oriented trees]}, Version v1.0.0, Zenodo, doi:10.5281/zenodo.22652747.
		
		\bibitem{doi: 10.1137/22M148046X}
		Jos\'{e} Aliste-Prieto, Anna de~Mier, Rosa Orellana, and Jos\'{e} Zamora.
		\newblock Marked graphs and the chromatic symmetric function.
		\newblock {\em SIAM Journal on Discrete Mathematics}, 37(3): 1881--1919, 2023.
		
		\bibitem{campbell2024liftchromaticsymmetricfunctions}
		John~M. Campbell.
		\newblock A lift of chromatic symmetric functions to $\textsf{NSym}$, 2024,arXiv:2410.04669.
		
		\bibitem{duchamp: hal-00018540}
		G{\'e}rard Henry~Edmond Duchamp, Alexander Klyachko, Daniel Krob, and Jean-Yves Thibon.
		\newblock {Noncommutative symmetric functions III :  Deformations of Cauchy and convolution algebras}.
		\newblock {\em {Discrete Mathematics and Theoretical Computer Science}}, Vol. 1: 159--216, January 1997.
		\newblock [in ''Special Issue :  Lie Computations'', G. Jacob, V. Koseleff, Eds.].
		
		\bibitem{gelfand1994noncommutativesymmetricfunctions}
		Israel Gelfand, D.~Krob, Alain Lascoux, B.~Leclerc, V.~S. Retakh, and J.~Y. Thibon.
		\newblock Noncommutative symmetric functions, 1994.
		
		\bibitem{heil2018algorithmcomparingchromaticsymmetric}
		Sam Heil and Caleb Ji.
		\newblock On an algorithm for comparing the chromatic symmetric functions of trees, 2018.
		
		\bibitem{AIHPD_2019__6_3_357_0}
		Martin Loebl and Jean-S\'ebastien Sereni.
		\newblock Isomorphism of weighted trees and {Stanley's} isomorphism conjecture for caterpillars.
		\newblock {\em Annales de l{\textquoteright}Institut Henri Poincar\'e D}, 6(3): 357--384, 2019.
		
		\bibitem{10.1093/oso/9780198534891.001.0001}
		I~G Macdonald.
		\newblock {\em Symmetric Functions and Hall Polynomials}.
		\newblock Oxford University Press, 03 1995.
		
		\bibitem{MARTIN2008237}
		Jeremy~L. Martin, Matthew Morin, and Jennifer~D. Wagner.
		\newblock On distinguishing trees by their chromatic symmetric functions.
		\newblock {\em Journal of Combinatorial Theory, Series A}, 115(2): 237--253, 2008.
		
		\bibitem{STANLEY1995166}
		R.P. Stanley.
		\newblock A symmetric function generalization of the chromatic polynomial of a graph.
		\newblock {\em Advances in Mathematics}, 111(1): 166--194, 1995.
		
		\bibitem{wang2024classtreesdeterminedchromatic}
		Yuzhenni Wang, Xingxing Yu, and Xiao-Dong Zhang.
		\newblock A class of trees determined by their chromatic symmetric functions, 2024.
		
	\end{thebibliography}

\end{document}